\documentclass[a4paper,10pt]{article}
\usepackage{amsmath}
\usepackage{amssymb,amsfonts,amsthm}
\usepackage[utf8]{inputenc} %acentos
\usepackage[english]{babel}

\usepackage{anysize}
\usepackage{enumerate}
\usepackage{tikz}
\usepackage{subfig}
\usepackage[all]{xy}
\usepackage{mathrsfs}
\usepackage{verbatim}
\usepackage{cite}
\definecolor{green}{rgb}{0,0.8,0.5}

\usepackage{hyperref}
\hypersetup{colorlinks=true,linkcolor=blue,citecolor=red}
\usepackage[width=0.8\textwidth]{caption}

\renewenvironment{abstract}{\small\quotation\noindent
 {\bfseries \abstractname .}}{\endquotation \par}

\newenvironment{trivlistparm}[2]{\trivlist
\item[{#1 #2}]}{\endtrivlist}

\newenvironment{prooftext}[1]{\trivlistparm{\bfseries}{#1}}{\Qed\endtrivlistparm}
\newenvironment{prova}{\trivlistparm{\bfseries}{Proof.}}{\Qed\endtrivlistparm}

\catcode`\@=11

\def\resetthefootnote{\renewcommand{\thefootnote}{\@arabic\c@footnote} }
\def\@principiremex#1{\trivlist
 \item[\hskip \labelsep{\bfseries #1\ \thethm.}]\ignorespaces}
\def\opar@principiremex#1[#2]{\trivlist
 \item[\hskip \labelsep{\bfseries #1\ \thethm\ (#2).}]\ignorespaces}

\newcommand{\newTHEOremrom}[2]{\newenvironment{#1}{\refstepcounter{thm}\@ifnextchar[{\opar@principiremex{#2}}
{\@principiremex{#2}}}{\qedB\endtrivlist}} \catcode`\@=12
\DeclareMathSymbol{\square}{\mathord}{AMSa}{"03}
\newcommand{\qedB}{\nopagebreak\hspace*{\fill}$\square$\par}
\newcommand{\Qed}{\nopagebreak\hspace*{\fill}{\vrule width6pt height6pt depth0pt}\par}

\newTHEOremrom{defi} {Definition}
\newTHEOremrom {rem} {Remark}
\newTHEOremrom {ex} {Example}

\newcommand{\refc}[1]{\mbox{$(\ref{#1})$}}
\newcommand{\secc}[1]{Section~\ref{#1}}
\newcommand{\teoc}[1]{Theorem~\ref{#1}}

\newcommand{\lemc}[1]{Lemma~\ref{#1}}
\newcommand{\defic}[1]{Definition~\ref{#1}}

\newcommand{\figc}[1]{Figure~\ref{#1}}

\renewcommand{\geq}{\geqslant}
\renewcommand{\epsilon}{\varepsilon}
\renewcommand{\leq}{\leqslant}
\newcommand{\F}{\mathscr{F}}
\newcommand{\LL}{\mathscr{L}}
\newcommand{\CC}{\mathscr{B}}
\newcommand{\PP}{\mathcal{P}}
\newcommand{\DD}{\mathscr{D}}
\newtheorem{thm}{Theorem}[section]
\newtheorem{thmx}{Theorem}
\newtheorem{prop}[thm]{Proposition}

\newtheorem{lema}[thm]{Lemma}

\newcommand{\abs}[1]{\left|#1\right|}
\newcommand{\R}{\mathbb{R}}

\newcommand{\Z}{\mathbb{Z}}
\newcommand{\N}{\mathbb{N}}

\newcommand{\C}{\mathbb{C}}
\newcommand{\PA}{\mathscr{P}}
\newcommand{\PI}{\mathcal{I}}

\newcommand{\RP}{\mathbb{RP}}

\newcommand{\out}{\Pi}%{\mathscr E}
\newcommand{\fun}{\mathscr G}

\newcommand{\sist}[2]{
  \left\{\!
   \begin{array}{l}
    \dot x=#1 \\[2pt] \dot y=#2
   \end{array}
  \right.
}

\title{\textbf{On the upper bound of the criticality of potential systems at the outer boundary using the Roussarie-Ecalle compensator}
\footnotetext{{\it E-mail address}:  rojas@ugr.es.}
}

\author{David Rojas\\[10pt]
{\small \textsl{Departamento de Matem\'atica Aplicada, Universidad de Granada, 18071 Granada, Spain.}}\\
{\small \textsl{Departament de Matem\`atiques, Universitat Aut\`onoma de Barcelona, 08193 Bellaterra, Spain.}}\\[5pt]
}

\date{}

\begin{document}

\maketitle

\begin{abstract}
This paper is concerned with the study of the criticality of families of planar centers. More precisely, we study sufficient conditions to bound the number of critical periodic orbits that bifurcate from the outer boundary of the period annulus of potential centers. In the recent years, the new approach of embedding the derivative of the period function into a collection of functions that form a Chebyshev system near the outer boundary has shown to be fruitful in this issue. In this work, we tackle with a remaining case that was not taken into account in the previous studies in which the Roussarie-Ecalle compensator plays an essential role. The theoretical results we develop are applied to study the bifurcation diagram of the period function of two different families of centers: the power-like family $\ddot x=x^p-x^q$, $p,q\in\R$ with $p>q$; and the family of dehomogenized Loud's centers.
\end{abstract}

\section{Introduction}\label{sec:intro}

The present paper deals with the bifurcation of critical periodic orbits of families of planar potential centers. Let $\Lambda$ be an open set of $\R^d$, $d\geq 1$, $\mu\in\Lambda$ be a parameter and consider a continuous family of analytic functions $V_{\mu}=V_{\mu}(x)$ defined in an open interval $I_{\mu}\subset \R$ that contains $x=0$. It is well known that the system
\begin{equation}\label{sys}
\dot x = y, \ \dot y = -V_{\mu}'(x)
\end{equation}
has a non-degenerate center at the origin for each $\mu\in\Lambda$ if $V_{\mu}(0)=V_{\mu}'(0)=0$ and $V_{\mu}''(0)>0$. That is, the origin has a punctured neighborhood entirely foliated by closed orbits surrounding it. The largest neighborhood with this property is called the period annulus of the center and we will denote it by $\PA_{\mu}$. After embedding $\PA_{\mu}$ into $\RP^2$, the boundary of the period annulus has two connected components: the center itself (which is called the inner boundary) and the outer boundary defined by $\Pi_{\mu}\!:=\partial\PA_{\mu}\setminus\{(0,0)\}$. The periodic orbits are inside the energy levels of the Hamiltonian $H(x,y;\mu)=\frac{y^2}{2}+V_{\mu}(x)$. We have then $H(\PA_{\mu})=(0,h_0(\lambda))$ with $h_0(\mu)\in\R^+\cup\{+\infty\}$. The inner boundary is inside the energy level $h=0$ for all $\mu\in\Lambda$ and we shall say that $h_0(\mu)$ is the energy level of the outer boundary $\Pi_{\mu}$. The minimal period $T_{\mu}(h)$ of the periodic orbit $\gamma_{h,\mu}$ inside the energy level $\{H(x,y;\mu)=h\}$ is given by the Abelian integral
\[
T_{\mu}(h)=\int_{\gamma_{h,\mu}}\frac{dx}{y}.
\]
This function is analytic on $(0,h_0(\mu))$ for each $\lambda\in\Lambda$ and it can be extended analytically to $h=0$ since the center is non-degenerate. Its derivative $T_{\mu}'(h)$ is also given by an Abelian integral and this paper is concerned with its zeros near the energy level $h=h_0(\mu)$, which correspond to critical periodic orbits near $\Pi_{\mu}$. More concretely, for a fixed $\hat\mu\in\Lambda$, we study the number of critical periodic orbits of system~\eqref{sys} that can emerge or disappear from $\Pi_{\hat \mu}$ as we move slightly the parameter $\mu\approx\hat\mu$. This number is called the criticality of the outer boundary.

\begin{defi}\label{defi2}
Consider a continuous family $\{X_{\mu}\}_{\mu\in\Lambda}$ of planar analytic vector fields with a center and fix some $\hat\mu\in\Lambda$. Suppose that the outer boundary of the period annulus varies continuously at $\hat\mu\in\Lambda$, meaning that for any $\varepsilon>0$ there exists $\delta>0$ such that $d_{H}(\out_{\mu},\out_{\hat\mu})\leqslant\varepsilon$ for all $\mu\in\Lambda$ with $\|\mu-\hat\mu\|\leqslant\delta$. Then, setting 
\[
 N(\delta,\varepsilon)\!:=\sup\left\{\text{\# critical periodic orbits $\gamma$ of $X_{\mu}$ in $\PA_{\mu}$ with $d_{H}(\gamma,\out_{\hat\mu})\leqslant\varepsilon$ and $\|\mu-\hat\mu\|\leqslant\delta$}\right\},
\]
the \emph{criticality} of $(\out_{\hat\mu},X_{\hat\mu})$ with respect to the deformation $X_{\mu}$ is $\mathrm{Crit}\bigl((\out_{\hat\mu},X_{\hat\mu}),X_{\mu}\bigr)\!:=\inf_{\delta,\varepsilon}N(\delta,\varepsilon).$
\end{defi}

In the previous definition $d_H$ stands for the Hausdorff distance between compact sets of $\RP^2$. The criticality of $(\out_{\hat\mu},X_{\hat\mu})$ may be infinite but in the case it is finite it gives the maximal number of critical periodic orbits of system $X_{\mu}$ that tend to the outer boundary $\Pi_{\hat\mu}$ in the Hausdorff sense as the parameter $\mu$ approach $\hat\mu$. We stress the requirement of the assumption that the period annulus varies continuously, which ensures that the changes on the geometry of $\PA_{\mu}$ do not occur abruptly as we vary the parameters of the system (see \cite{ManVil06} for details). On account of the previous definition, we say that a parameter $\hat\mu\in\Lambda$ is a local regular value of the period function at the outer boundary of the period annulus if $\mathrm{Crit}\bigl((\out_{\hat\mu},X_{\hat\mu}),X_{\mu}\bigr)=0.$ Otherwise we say that it is a local bifurcation value at the outer boundary.

The study of critical periodic orbits is analogous to the study of limit cycles, the objects of main concern of the Hilbert's 16th problem (see for instance~\cite{Blows-Lloyd,Pugh,Roussarie,Ye} and references there in). Questions related to the behavior of the period function have been extensively studied by a large number of authors. Let us quote for instance the problems of isochronicity (see \cite{CMV,loud,Pavao1}), monotonicity (see \cite{ChiconeI,ChiconeII,schaaf}) or bifurcation of critical periodic orbits (see \cite{CJ,RT,smoller}).

In the collection of works~\cite{ManRojVil2015,ManRojVil2016-JDDE} tools that enable to bound the criticality at the outer boundary for the family of potential systems~\eqref{sys} were developed. These tools, and the ones we present here, allow to tackle the bifurcation problem in the following two situations: either $h_0(\mu)=+\infty$ or $h_0(\mu)<+\infty$ for all $\mu\approx\hat\mu$. (The case in which in any neighborhood of $\hat\mu$ there are $\mu_1$ and $\mu_2$ with $h_0(\mu_1)=+\infty$ and $h_0(\mu_2)<+\infty$ is not considered.) For each one of these two situations, sufficient conditions in order that $\mathrm{Crit}\bigl((\out_{\hat\mu},X_{\hat\mu}),X_{\mu}\bigr)\leq n$ for $n\in\N\cup\{0\}$ were given in~\cite[Theorem A and Theorem B]{ManRojVil2016-JDDE}. The idea in both cases was to find a collection of functions $\phi_{\mu}^i(h)$, $i=1,2,\ldots,n$, verifying that there exist $\delta,\varepsilon>0$ such that $(\phi_{\mu}^1,\phi_{\mu}^2,\ldots,\phi_{\mu}^n,T_{\mu}')$ form an Extended Complete Chebyshev system (ECT-system for short, see \defic{defi_ECT}) on the interval $(h_0(\mu)-\varepsilon,h_0(\mu))$ if $\|\mu-\hat\mu\|<\delta$. In particular this implies that $T_{\mu}'(h)$ has at most $n$ zeros in $(h_0(\mu)-\varepsilon,h_0(\mu)),$ counted with multiplicities, uniformly for all parameters $\mu\approx\hat\mu$ and so
$\mathrm{Crit}\bigl((\out_{\hat\mu},X_{\hat\mu}),X_{\mu}\bigr)\leqslant n.$ The problem consisted then to guarantee that the Wronskian (see \defic{wronskian}) of $(\phi_{\mu}^1,\phi_{\mu}^2,\ldots,\phi_{\mu}^n,T_{\mu}')$ is different from zero for all $h\approx h_0(\mu)$ and $\mu\approx\hat\mu.$ 
In the present paper we extend the results~\cite[Theorems~A and~B]{ManRojVil2016-JDDE} by considering a remaining case which at that moment the techniques did not cover (see Theorems~\ref{thm:criticality_infinite} and~\ref{thm:criticality_finite}). To do so, additional tools to the ones in these works are developed in Section~\ref{sec:technical}. The treatment of the necessary conditions to bound the criticality in this limit cases are presented in Section~\ref{sec:criticality}. 

\begin{figure}
\centering
\subfloat[Previous bifurcation diagram\label{fig:diagram}]
{
\centering \includegraphics[scale=1]{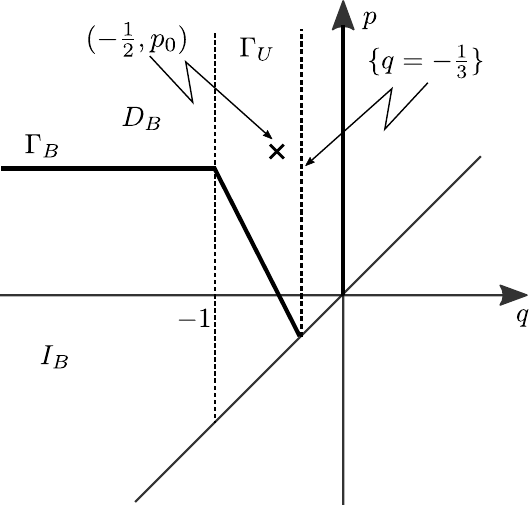}
}
\hspace{1.5cm}
\subfloat[New bifurcation diagram\label{fig:diagram-nou}]
{
\centering \includegraphics[scale=1]{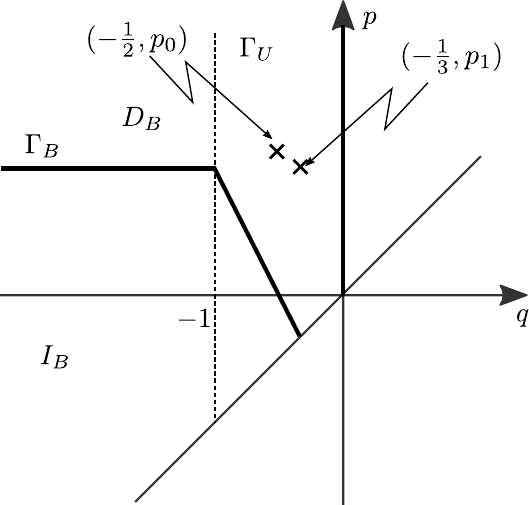}
}
\caption{On the left, bifurcation diagram of the period function of the family~\eqref{Xmu2} at the outer boundary of the period annulus according to \cite[Theorem~E]{ManRojVil2015} and \cite[Theorem~C]{ManRojVil2016-JDDE}. On the right, improvement of the bifurcation diagram according with Theorem~\ref{thm:familia}. In both figures, $\Gamma_B$ and $\Gamma_U$ stand, respectively, for the union of the solid and dotted lines.}
\end{figure}

As in the previous works, our testing ground is the two-parametric family of potential differential systems given by
\begin{equation}\label{Xmu2}
\begin{cases}
  \dot{x}=-y,\\
 \dot{y}=(x+1)^p-(x+1)^q,
 \end{cases}
\end{equation}
which has a non-degenerate center at the origin for all $\mu\!:=(q,p)$ varying in $\Lambda\!:=\{(q,p)\in\R^2:p>~q\}$. Note that, for each $\mu\in\Lambda,$  system~\eqref{Xmu2} is analytic on $\{(x,y)\in\R^2:x>-1\}.$ Our interest in this family began because of the results by Miyamoto and Yagasaki \cite{MiyYag13} about the monotonicity of the period function for $q=1$ and $p\in\N$. Later Yagasaki \cite{Yagasaki13} improved this result proving the monotonicity property remains if one consider any $p>1$ real. 
Following that, we performed a more exhaustive study of the period function of the family~\eqref{Xmu2} for all $\mu\in\Lambda$ in~\cite{ManRojVil2015,ManRojVil2016-JMAA,ManRojVil2016-JDDE}. To be more precise, in~\cite{ManRojVil2016-JMAA} we were concerned with the monotonicity of the period function, the criticality of the inner boundary and the criticality of the interior of the period annulus of its isochronous centers. In~\cite{ManRojVil2015,ManRojVil2016-JDDE} we studied the criticality of the outer boundary and it is precisely the result we obtained there for the family~\refc{Xmu2} the one that we improve here. In short, see \figc{fig:diagram}, we proved that $\mathrm{Crit}\bigl((\out_{\hat\mu},X_{\hat\mu}),X_{\mu}\bigr)=0$ if $\hat\mu\in\Lambda\setminus\{\Gamma_B\cup\Gamma_U\cup\{(-\tfrac{1}{2},p_0)\}$ with $p_0\approx 1.20175$, and $\mathrm{Crit}\bigl((\out_{\hat\mu},X_{\hat\mu}),X_{\mu}\bigr)\geqslant 1$ if $\hat\mu\in\Gamma_B.$ Moreover, we showed that the criticality is exactly one for parameters $\hat\mu\in\Gamma_B$ such that $\hat\mu=(0,\hat p)$ with $\hat p\in(0,+\infty)\setminus\{\tfrac{1}{2},1\}$, $\hat\mu=~(\hat q,1)$ with $\hat q<-3$ and $\hat\mu=(\hat q,-2\hat q-1)$ with $\hat q\in(-\tfrac{3}{5},-\tfrac{1}{3})\setminus\{-\tfrac{1}{2}\}$.
By applying the new tools in this paper we can go further, see Figure~\ref{fig:diagram-nou}, and prove the following result, where $f(p)\!:=(\tfrac{3+3p}{2})^{\frac{3+3p}{1+3p}}-3p-2$.

\begin{thmx}\label{thm:familia}
Let $\{X_{\mu}\}_{\mu\in\Lambda}$ be the family of potential vector fields in~\eqref{Xmu2} and consider the period function of the center at the origin. Then the following hold:
\begin{enumerate}[$(a)$]
\item If $\hat\mu=(-\tfrac{1}{3},\hat p)$ with $\hat p\in(-\tfrac{1}{3},+\infty)\setminus\{p_1\}$, where $p_1\approx 1.15685$ is the unique zero of $f$ on $(-\tfrac{1}{3},+\infty)$, then $\mathrm{Crit}\bigl((\Pi_{\hat\mu},X_{\hat\mu}),X_{\mu}\bigr)=0$.
\item If $\hat\mu=(0,\hat p)$ with $\hat p\in(0,+\infty)\setminus\{1\}$ then $\mathrm{Crit}\bigl((\Pi_{\hat\mu},X_{\hat\mu}),X_{\mu}\bigr)=1$.
\end{enumerate}
\end{thmx}

\begin{figure}
\centering
\includegraphics[scale=1]{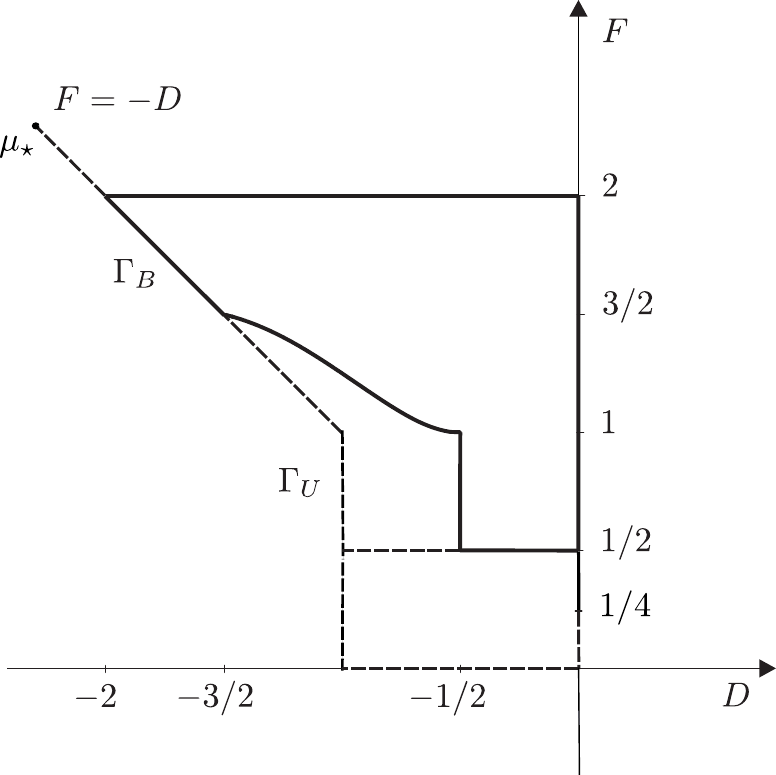}
\caption{\label{fig:diagrama-loud}Bifurcation diagram of the period function at the polycycle according to~\cite{ManVil08,MMV2,MMV4,MarMarSaaVil2015,MarVil06,Vil07}, where $\mu_{\star}=(-F_{\star},F_{\star})$ with $F_{\star}\approx 2.34$. The union of the bold curves correspond to the set of bifurcation parameters at the outer boundary. The union of dotted straight lines correspond to the set of unspecified parameters.}
\end{figure}

We finish this section stating a second application of the techniques developed in this work. The results obtained in the series of papers \cite{ManVil08,MMV2,MMV4,MarMarSaaVil2015,MarVil06,Vil07} by Ma\~nosas, Marde{\v{s}}i{\'c}, Mar\'in, Saavedra and Villadelprat deal also with the bifurcation of critical periodic orbits from the outer boundary of the period annulus. Their testing ground is the family of demohogenized Loud's centers 
\begin{equation}\label{loud}
\sist{-y+xy,}{x+Dx^2+Fy^2,}
\end{equation}
where $\mu\!:=(D,F)\in\R^2$. In this collection of works, the bifurcation diagram of the period function at the polycycle has been studied (see Figure~\ref{fig:diagrama-loud}). We refer the reader to the previous works for complete information and to~\cite{RojVil2018} for a summary of the latest results. We point out that the Loud's family can be brought to potential form by means of an explicit coordinate transformation, see \cite[Lemma 2.2]{Vil07}, and hence it is susceptible to be studied with our methods. In this paper we prove the following result.

\begin{thmx}\label{thm:loud}
Let $\{X_{\mu}\}_{\mu\in\Lambda}$ be the family of potential vector fields in~\eqref{loud} and consider the period function of the center at the origin. If $\hat\mu=(\hat D, 2)$ with $\hat D\in(-2,0)\setminus\{-1/2\}$ then $\mathrm{Crit}\bigl((\Pi_{\hat\mu},X_{\hat\mu}),X_{\mu}\bigr)=1$.
\end{thmx}

\section{Technical results}\label{sec:technical}

In this section we develop the technical tools that we will use in Section~\ref{sec:criticality} to prove the results concerning the criticality. Let $a\in\R^+\cup\{+\infty\}$ and consider the integral operator
\[
\F: C^{\omega}[0,a)\rightarrow C^{\omega}[0,a)
\]
defined by
\begin{equation}\label{F}
\F[f](x)\!:=\int_0^{\frac{\pi}{2}}f(x\sin\theta)d\theta.
\end{equation}
Here, and it what follows, $C^{\omega}[0,a)$ stands for the set of analytic functions on $(0,a)$ that can be extended analytically to $x=0$. In \cite{ManRojVil2016-JDDE} it was shown that the derivative of the period function of system~\eqref{sys} is related with this operator by means of the equality
\begin{equation}\label{formula}
\sqrt{2}h^2T_{\mu}'(h^2)=\F[f_{\mu}](h),
\end{equation}
where $f_{\mu}(x)=x(g_{\mu}^{-1})''(x)-x(g_{\mu}^{-1})''(-x)$ and $g_{\mu}(x)\!:=\text{sgn}(x)\sqrt{V_{\mu}(x)}$. Note that $f_{\mu}$ is in the class $C^{\omega}[0,h_0(\mu))$. 
Let us recall the notions of Chebyshev system and Wronskian, that will be useful for our purposes.   

\begin{defi}\label{defi_ECT}
Let $f_0,f_1,\dots f_{n-1}$ be analytic functions on an open real interval $I$. The ordered set
$(f_0,f_1,\dots f_{n-1})$ is an \emph{extended complete Chebyshev system} (for short, a ECT-system) on $I$ if, for all $k=1,2,\dots n$, any nontrivial linear combination
\[
\alpha_0f_0(x)+\alpha_1f_1(x)+\cdots+\alpha_{k-1}f_{k-1}(x)
\]
has at most $k-1$ isolated zeros on $I$ counted with multiplicities.
(In these abbreviations, ``T'' stands for Tchebycheff, which in some sources is the transcription of the Russian name Chebyshev).
\end{defi}

\begin{defi}\label{wronskian}
Let $f_0,f_1,\dots,f_{k-1}$ be analytic functions on an open interval $I$ of $\R$. Then
\[
W[f_0,f_1,\dots,f_{k-1}](x)\!:=\text{det}\left(f_j^{(i)}(x)\right)_{0\leq i,j\leq k-1}=
\left|\begin{array}{ccc}
f_0(x) & \cdots & f_{k-1}(x)\\
f_0'(x) & \cdots & f_{k-1}'(x)\\
		&	\vdots & \\
f_0^{(k-1)}(x) & \cdots & f_{k-1}^{(k-1)}(x)
\end{array}\right|
\]
is the \emph{Wronskian} of $(f_0,f_1,\dots,f_{k-1})$ at $x\in I$.
\end{defi}

The previous two notions are closely related by the following result (see \cite{Karlin}).

\begin{lema}\label{lema:ECT-Wronskia}
$(f_0,f_1,\dots,f_{n-1})$ is an ECT-system on $I$ if and only if, for each $k=1,2,\dots,n$,
\[
W[f_0,f_1,\dots,f_{k-1}](x)\neq 0\ \ \text{for all }x\in I.
\]
\end{lema}

With equality~\eqref{formula} in mind, the aim is to complete $\F[f_{\mu}]$ with a collection of analytic functions $\phi_{\mu}^1,\dots,\phi_{\mu}^n$ in order that $(\phi_{\mu}^1,\dots,\phi_{\mu}^n,\F[f_{\mu}])$ form an ECT-system on $(a-\epsilon,a)$ for some $\epsilon>0$, for all $\mu\approx\hat\mu$. Notice that, to obtain the desired upper bounds on the criticality, the crucial point is to guarantee the uniformity with respect to the parameters of the system. In the works~\cite{ManRojVil2015,ManRojVil2016-JDDE} sufficient conditions in terms of $f_{\mu}$ in order that $\F[f_{\mu}]$ can be embedded into an ECT-system were given. These conditions were formulated using the notions that we introduce next.
 
\begin{defi}\label{def:quantifiable}
Let $f\in C^{\omega}(a,b)$. We say that $f$ is \emph{quantifiable at $b$ by $\alpha$ with limit $\ell$} in case that:
\begin{enumerate}[$(i)$]
\item If $b\in\R$, then $\lim_{x\rightarrow b^-} f(x)(b-x)^{\alpha}=\ell$ and $\ell\neq 0$.
\item If $b=+\infty$, then $\lim_{x\rightarrow+\infty} \frac{f(x)}{x^{\alpha}}=\ell$ and $\ell\neq 0$.
\end{enumerate}
We call $\alpha$ the \emph{quantifier} of $f$ at $b$. We shall use the analogous definition at $a$.
\end{defi}

\begin{defi}\label{def:cont_quantifiable}
Let $\{f_{\mu}\}_{\mu\in\Lambda}$ be a continuous family of analytic
functions on~$(a(\mu),b(\mu))$, meaning that the map $(x,\mu)\longmapsto f_{\mu}(x)$ is continuous on $\{(x,\mu)\in\R\times\Lambda: x\in (a(\mu),b(\mu))\}$. Assume that $b$ is either a continuous function from $\Lambda$ to $\R$ or $b(\mu)= +\infty$ for all $\mu\in\Lambda$. Given $\hat\mu\in\Lambda$
we shall say that $\{f_{\mu}\}_{\mu\in\Lambda}$ is
\emph{continuously quantifiable} in $\hat\mu$ at $b(\mu)$ by
$\alpha(\mu)$ with limit $\ell(\hat\mu)$ if there exists an open neighborhood~$U$ of $\hat\mu$ such that $f_{\mu}$ is
quantifiable at $b(\mu)$ by $\alpha(\mu)$ with limit $\ell(\mu)$ for all $\mu\in U$ and, moreover,
\begin{enumerate}[$(i)$]
\item In case that $b(\hat\mu)<+\infty$, then $\lim_{(x,\mu)\rightarrow (b(\hat\mu),\hat\mu)} f_{\mu}(x)(b(\mu)-x)^{\alpha(\mu)}=\ell(\hat\mu)$ and $\ell(\hat\mu)\neq 0$.
\item In case that $b(\hat\mu)=+\infty$, then $\lim_{(x,\mu)\rightarrow (+\infty,\hat\mu)} \frac{f_{\mu}(x)}{x^{\alpha(\mu)}}=\ell(\hat\mu)$ and $\ell(\hat\mu)\neq 0$.
\end{enumerate}
For the sake of shortness, in the first case we shall write $f_{\mu}(x)\sim_{b(\mu)}\ell(\mu)(b(\mu)-x)^{-\alpha(\mu)}$ at $\hat\mu$, and in the second case $f_{\mu}(x)\sim_{+\infty}\ell(\mu)x^{\alpha(\mu)}$ at $\hat\mu$.
We shall use the analogous definition for the left endpoint $a(\mu)$. 
\end{defi}

We point out that the map  $\alpha\!:U\longrightarrow \R$ that appears in the previous definition must be continuous at~$\hat\mu$ (see~\cite[Remark 2.6]{ManRojVil2016-JDDE}). 
Given $\nu_1,\nu_2,\ldots,\nu_n\in\R$, consider the linear ordinary differential operator 
\[
\LL_{\boldsymbol\nu_n}\!:\mathscr C^{\omega}(0,+\infty)\longrightarrow \mathscr C^{\omega}(0,+\infty)
\]
given by
\begin{equation}\label{eq:op_Wronskia}
\LL_{\boldsymbol\nu_n}[f](x)\!:=\frac{W\left[x^{\nu_1},x^{\nu_2},\ldots,x^{\nu_n},f(x)\right]}{x^{\sum_{i=1}^n(\nu_i-i)}}.
\end{equation}
Here, and in what follows, for the sake of shortness we use the notation $\boldsymbol\nu_n=(\nu_1,\dots,\nu_n)$. Furthermore we define $\LL_{\boldsymbol\nu_0}=id$ in order that the statements of the next results contemplate the case $n=0$ as well.

From now one and till the end of this section let us assume that $a=+\infty$ and so $f_{\mu}\in C^{\omega}[0,+\infty)$. In~\cite{ManRojVil2015,ManRojVil2016-JDDE} sufficient conditions on the family $\{\LL_{\boldsymbol\nu_n}[f_{\mu}]\}_{\mu\in\Lambda}$ in order that the family $\{(\LL_{\boldsymbol\nu_n}\circ\F)[f_{\mu}]\}_{\mu\in\Lambda}$ is continuously quantifiable at infinity were given. One of the major requirements was that the continuous family $\{\LL_{\boldsymbol\nu_n}[f_{\mu}]\}_{\mu\in\Lambda}$ satisfies $\LL_{\boldsymbol\nu_n}[f_{\mu}](x)\sim_{\infty}\ell(\mu)x^{\alpha(\mu)}$ at $\hat\mu$ with $\alpha(\hat\mu)\neq -1+2m$, $m\in\N\cup\{0\}$. This fact, among other specific assumptions, allowed to bound the number of zeros of $(\LL_{\boldsymbol\nu_n}\circ\F)[f_{\mu}](x)$ near infinity uniformly on the parameters $\mu\approx\hat\mu$. Consequently, criticality results were obtained on account of Lemma~\ref{lema:ECT-Wronskia} (see~\cite[Proposition 2.16 and Theorem A]{ManRojVil2016-JDDE}). 
In the present paper we aim to extend these results to the case $\alpha(\hat\mu)=-1+2m$, $m\in\N\cup\{0\}$. This situation presents a much more complicated behavior than the general one. 

For the sake of simplicity in the forthcoming illustration, let restrict ourselves to the case $n=0$ and $m=0$, and consider the continuous family of analytic functions $f_{\mu}(x)=x^{\alpha(\mu)}$. 
In \cite[Theorem~2.13]{ManRojVil2015} it was proved that if $\alpha(\hat\mu)>-1$ then $\F[f_{\mu}](x)\sim_{\infty}x^{\alpha(\mu)}$ at $\hat\mu$. In addition, using \cite[Theorem~2.17]{ManRojVil2015} one can easily show that if $\alpha(\hat\mu)\in(-2,-1)$ then $\F[f_{\mu}](x)\sim_{\infty}\frac{1}{x}$ at $\hat\mu$. These two results are uniformly with respect to the parameters. However, if one fix $\mu=\hat\mu$ such that $\alpha(\hat\mu)=-1$ then, see~\cite[Proposition 2.5]{ManRojVil2015},
\[
\lim_{x\rightarrow+\infty}\frac{x}{\log(x)}\F[f_{\hat\mu}](x)=1.
\]
This situation can not be covered uniformly on the parameters using the notions introduced in Definition~\ref{def:cont_quantifiable} because of the appearance of a logarithm term. In order to deal with this, we shall use the so-called Roussarie-Ecalle compensator (see~\cite{Roussarie}).

\begin{defi}\label{defi:compensador}
The function defined for all $x>0$ and $\alpha\in\R$ by means of
\[
\omega(x,\alpha)\!:= \begin{cases}
\frac{x^{\alpha+1}-1}{\alpha+1} & \text{ if }\alpha\neq -1,\\
\log x & \text{ if }\alpha=-1.
\end{cases}
\]
is called the \emph{Roussarie-Ecalle compensator}.
\end{defi}

Here we centered the singular value of the parameter at $\alpha=-1$ instead of the classical version at $\alpha=0$ for the sake of simplicity through this paper. A useful property that was proved in~\cite{Roussarie} and we shall use here is that
\begin{equation}\label{limit_omega}
\lim_{(x,\alpha)\rightarrow (+\infty,-1)}\omega(x,\alpha)=+\infty.
\end{equation}
With the purpose to state the main result of this section properly, we denote
\begin{equation}\label{G}
\fun(\alpha)\!:=\frac{\sqrt{\pi}}{2}\frac{\Gamma\left(\frac{1+\alpha}{2}\right)}{\Gamma\left(1+\frac{\alpha}{2}\right)}
\end{equation}
where $\Gamma$ is the Gamma function and
\begin{equation}\label{K}
\mathscr K(\alpha)\!:=\begin{cases}
\fun(\alpha)-\frac{1}{1+\alpha} & \text{ if }\alpha\neq -1,\\
\log(2) & \text{ if }\alpha=-1,
\end{cases}
\end{equation}
which is a continuous positive function on $(-2,0)$ (see Lemma~\ref{lema:funcions}). For the sake of brevity, we also define
\begin{equation}\label{Omega}
\Omega(x,\alpha)\!:=(\alpha+1)\fun(\alpha)\omega(x,\alpha),
\end{equation}
which is also a continuous positive function for all $x>1$ and $\alpha>-2$ (see Lemma~\ref{lema:funcions}).
Following the notation in Definition~\ref{def:cont_quantifiable}, we shall write $f_{\mu}(x)\sim_{+\infty}\ell(\mu)\Omega(x,\alpha(\mu))$ at $\hat\mu$ if 
\[
\lim_{(x,\mu)\rightarrow(+\infty,\hat\mu)}\frac{1}{\Omega(x,\alpha(\mu))}f_{\mu}(x)=\ell(\hat\mu)\neq 0.
\]

\begin{defi}
Let $f\in C^{\omega}[0,+\infty)$. Then, for each $n\in\N$, we call
\[
M_n[f]\!:=\int_0^{+\infty}x^{2n-2}f(x)dx
\]
the $n$-th momentum of $f$, whenever it is well defined.
\end{defi}

Let us finally denote $b_0(\alpha)\equiv 1$ and $b_m(\alpha)\!:=\prod_{i=1}^m \frac{\alpha+2i}{\alpha+2i-1}$ for $m\geq 1$. In the next statement, as well as in the following one, the condition $M_1\equiv M_2\equiv\cdots\equiv M_m\equiv 0$ is void in case that $m=0$.

\begin{thmx}\label{thm:main}
Let $\Lambda$ be an open subset of $\R^d$ and consider a continuous family $\{f_{\mu}\}_{\mu\in\Lambda}$ of analytic functions on $[0,+\infty)$ such that $f_{\mu}(x)\sim_{\infty} a(\mu)x^{\alpha(\mu)}$ at $\hat\mu\in\Lambda$. Let us assume that there exists $m\in\N\cup\{0\}$ such that $\alpha(\hat\mu)+2m=-1$. If $M_1[f_{\mu}]\equiv M_2[f_{\mu}]\equiv\cdots\equiv M_m[f_{\mu}]\equiv 0$ then
\[
x^{2m+1}\F[f_{\mu}](x)\sim_{\infty} a(\mu)b_m(\alpha(\mu))\Omega(x,\alpha(\mu)+2m) \text{ at }\hat\mu.
\]
\end{thmx}

We point out at this point that Theorem~\ref{thm:main} covers the uniformity on the parameters that lacked in the above illustration example. Indeed, let us consider again $f_{\mu}(x)=x^{\alpha(\mu)}$ and let us fix $\hat\mu\in\Lambda$ such that $\alpha(\hat\mu)=-1$. According with Theorem~\ref{thm:main},
\[
\F[f_{\mu}](x)\sim_{\infty} \frac{1}{x}\Omega(x,\alpha(\mu)) \text{ at }\hat\mu.
\]
We invoke equality~\eqref{Omega} and Definition~\ref{defi:compensador} to show that
\[
\dfrac{1}{x}\Omega(x,\alpha(\mu))=
\begin{cases}
\dfrac{x^{\alpha(\mu)+1}-1}{x(\alpha(\mu)+1)}\sim_{+\infty}x^{\alpha(\mu)}  & \text{ if }\alpha(\mu)>-1,\\
\dfrac{\log x}{x}  & \text{ if }\alpha(\mu)=-1,\\
\dfrac{x^{\alpha(\mu)+1}-1}{x(\alpha(\mu)+1)}\sim_{+\infty} \dfrac{1}{x}  & \text{ if }-2<\alpha(\mu)<-1.\\
\end{cases}
\]

We finish this section stating a general version of Theorem~\ref{thm:main}. This generalization gives necessary conditions on the family $\{\LL_{\boldsymbol\nu_n}[f_{\mu}]\}_{\mu\in\Lambda}$ that enables to quantify $(\LL_{\boldsymbol\nu_n}\!\circ\F)[f_{\mu}]$ at $x=+\infty.$ In the following statement $\nu_1,\nu_2,\dots,\nu_n$ are not real numbers any more but continuous functions on $\Lambda$. For shortness, we keep using the notation $\boldsymbol\nu_n(\mu)=(\nu_1(\mu),\dots,\nu_n(\mu))$. Again, the condition of the existence of continuous function $\nu_1,\nu_2,\dots,\nu_n$ is void in case that $n=0$.

\begin{thmx}\label{thm:main-gen}
Let $\Lambda$ be an open subset of $\R^d$ and consider a continuous family $\{f_{\mu}\}_{\mu\in\Lambda}$ of analytic functions on $[0,+\infty)$. Assume that, in a neighborhood of some fixed $\hat\mu\in\Lambda$, there exist $n\geq 0$ continuous functions $\nu_1,\nu_2,\dots,\nu_n$, with $\nu_1(\hat\mu),\dots,\nu_n(\hat\mu)$ pairwise distinct, and $m\in\N\cup\{0\}$ such that 
$\LL_{\boldsymbol\nu_n}[f_{\mu}](x)\sim_{\infty} a(\mu)x^{\alpha(\mu)}$
at $\hat\mu$ with $\alpha(\hat\mu)+2m=-1$.  If $M_1[\LL_{\boldsymbol\nu_n}[f_{\mu}]]\equiv M_2[\LL_{\boldsymbol\nu_n}[f_{\mu}]]\equiv\cdots\equiv M_m[\LL_{\boldsymbol\nu_n}[f_{\mu}]]\equiv 0$ then
\[
x^{2m+1}(\LL_{\boldsymbol\nu_n}\circ\F)[f_{\mu}](x)\sim_{\infty}a(\mu)b_m(\alpha(\mu))\Omega(x,\alpha(\mu)+2m)
\text{ at }\hat\mu.
\]
\end{thmx}

\subsection{Proof of Theorem~\ref{thm:main}}
This section is devoted to the proof of Theorem~\ref{thm:main}. To this end, we first prove a collection of technical lemmas that will be useful for this purpose.

\begin{lema}\label{lema:funcions}
The following holds:
\begin{enumerate}[$(a)$]
\item For every fixed $x\in(1,+\infty)$ the function $\alpha\mapsto \Omega(x,\alpha)=(1+\alpha)\fun(\alpha)\omega(x,\alpha)$ is positive and monotonous increasing for all $\alpha\in(-2,+\infty)$. Moreover, $\lim_{\alpha\rightarrow -1}\Omega(x,\alpha)=\log(x)$ for all $x>1$.
\item The function $\alpha\mapsto \mathscr K(\alpha)$ is continuous and positive for all $\alpha\in(-2,+\infty)$.
\end{enumerate}
\end{lema}
\begin{prova}
Let us fix $x\in(1,+\infty)$. In order to prove $(a)$ we compute the derivative of the compensator using the expression on Definition~\ref{defi:compensador}. We have
\begin{equation}\label{f1}
(\alpha+1)\frac{d}{d\alpha}\omega=((\alpha+1)\omega+1)\log x-\omega.
\end{equation}
(Here we omit the dependence of $\omega$ with respect to $x$ and $\alpha$ for the sake of brevity.) Multiplying this expression by $(\alpha+1)$ and deriving again,
\[
\frac{d}{d\alpha}\left( (\alpha+1)^2\frac{d}{d\alpha}\omega\right)= (\alpha+1)((\alpha+1)\omega+1)\log^2 x.
\]
Since $(\alpha+1)\omega+1>0$ for all $\alpha>-2$, the previous expression vanishes only at $\alpha=-1$. This, together with the expression in~\eqref{f1} implies that 
$(\alpha+1)^2\frac{d}{d\alpha}\omega\geq 0$ for all $\alpha\in\R$ and the equality holds if and only if $\alpha=-1$. Deriving the expression in Definition~\ref{defi:compensador} we have that
\[
\lim_{\alpha\rightarrow -1}\frac{d}{d\alpha}\omega=\frac{1}{2}\log^2 x\neq 0.
\]
This shows the monotonicity of the map $\alpha\mapsto \omega(x,\alpha)$ for all $\alpha\in(-2,+\infty)$. Moreover, $\omega(x,-2)=1-\frac{1}{x}>0$. This implies that $\omega(x,\alpha)>0$ in $(-2,+\infty)$. Finally, on account of the properties of the Gamma function, the map
\[
\alpha\mapsto (1+\alpha)\fun(\alpha)=\frac{\sqrt{\pi}\Gamma\left(\frac{3+\alpha}{2}\right)}{\Gamma\left(1+\frac{\alpha}{2}\right)}
\]
is a well defined positive monotonous increasing function in $(-2,+\infty)$. Therefore $(1+\alpha)\fun(\alpha)\omega(x,\alpha)$ is a positive monotonous increasing function on $(-2,+\infty)$. Moreover, the previous equality implies that
\begin{equation}\label{limite}
\lim_{\alpha\rightarrow -1}(1+\alpha)\fun(\alpha)=1.
\end{equation}
Then $\lim_{\alpha\rightarrow -1}\Omega(x,\alpha)=\log(x)$ on account of Definition~\ref{defi:compensador}. This proves~$(a)$.

Let us now prove $(b)$. In order to see that $\mathscr K(\alpha)=\frac{(1+\alpha)\fun(\alpha)-1}{1+\alpha}>0$ we show that
\[
\frac{\sqrt{\pi}\Gamma\left(\frac{3+\alpha}{2}\right)}{\Gamma\left(1+\frac{\alpha}{2}\right)}>1 \text{ if }\alpha\in(-1,+\infty) \text{ and }0<\frac{\sqrt{\pi}\Gamma\left(\frac{3+\alpha}{2}\right)}{\Gamma\left(1+\frac{\alpha}{2}\right)}<1 \text{ if }\alpha\in(-2,-1).
\]
A direct computation shows that 
\begin{equation}\label{derivada-gamma}
\frac{d}{d\alpha}\!\!\left(\frac{\sqrt{\pi}\Gamma\left(\frac{3+\alpha}{2}\right)}{\Gamma\left(1+\frac{\alpha}{2}\right)}\right)=\frac{\sqrt{\pi}}{2}\frac{\Gamma\left(\frac{3+\alpha}{2}\right)}{\Gamma\left(1+\frac{\alpha}{2}\right)}\left(\psi\!\left(\tfrac{3+\alpha}{2}\right)-\psi(1+\tfrac{\alpha}{2})\right)
\end{equation}
where $\psi(x)=\Gamma'(x)/\Gamma(x)$ denotes the Digamma function (see~\cite[Section 6.3]{AbrSte92}). The function $\psi(x)$ is increasing for all $x>0$. Then, on account of the expression of the derivative, the function $\frac{\sqrt{\pi}\Gamma\left(\frac{3+\alpha}{2}\right)}{\Gamma\left(1+\frac{\alpha}{2}\right)}$ is increasing since it is positive for all $\alpha>-2$. The function $\mathscr K(\alpha)$ is positive then since 
\[
\lim_{\alpha\rightarrow -2^+}\frac{\sqrt{\pi}\Gamma\left(\frac{3+\alpha}{2}\right)}{\Gamma\left(1+\frac{\alpha}{2}\right)}=0 \text{ and }
\lim_{\alpha\rightarrow +\infty}\frac{\sqrt{\pi}\Gamma\left(\frac{3+\alpha}{2}\right)}{\Gamma\left(1+\frac{\alpha}{2}\right)}=+\infty.
\]
To show that $\lim_{\alpha\rightarrow -1}\bigl(\fun(\alpha)-\frac{1}{1+\alpha}\bigr)=\log(2)$ we invoke the limit~\eqref{limite} and, applying H\^opital's rule,
\[
\lim_{\alpha\rightarrow -1}\left(\fun(\alpha)-\frac{1}{1+\alpha}\right)=\lim_{\alpha\rightarrow -1}\sqrt{\pi}\frac{d}{d\alpha}\!\!\left(\frac{\Gamma\left(\frac{3+\alpha}{2}\right)}{\Gamma\left(1+\frac{\alpha}{2}\right)}\right)
\]
if this last limit exists. On account of equality~\eqref{derivada-gamma}, the result follows using that $\psi(1)-\psi(1/2)=2\log(2)$ (see~\cite[6.3.3]{AbrSte92}).
\end{prova}

\begin{lema}\label{lema:previ}
Let $\Lambda$ be an open subset of $\R^d$ and consider a continuous family $\{f_{\mu}\}_{\mu\in\Lambda}$ of analytic functions on $[0,+\infty)$. Let $\alpha(\mu)$ be a continuous function such that $\alpha(\hat\mu)=-1$ for some fixed $\hat\mu\in\Lambda$. Then for every $\epsilon>0$ and $M>0$ there exist $x_0>M$ and $\delta>0$ such that
\[
\abs{\frac{x}{\Omega(x,\alpha(\mu))}\int_0^{\arcsin(M/x)} f_{\mu}(x\sin\theta)d\theta}\leq \epsilon
\]
for all $x>x_0$ and $\|\mu-\hat\mu\|<\delta$.
\end{lema}
\begin{prova}
Let us fix $\epsilon,M>0$ and let $\delta>0$ such that $\alpha(\mu)\in(-2,0)$ for all $\|\mu-\hat\mu\|\leq \delta$. Since $\{f_{\mu}\}_{\mu\in\Lambda}$ is a continuous family of analytic functions on $[0,+\infty)$, $N=N(M,\delta)\!:=\max\{\abs{f_{\mu}(x)}: x\in[0,M], \|\mu-\hat\mu\|\leq \delta\}$ is well-defined. Then, for all $x>M$,
\[
\frac{x}{\Omega(x,\alpha(\mu))}\int_0^{\arcsin(M/x)}\abs{f_{\mu}(x\sin\theta)}d\theta \leq \frac{Nx\arcsin(M/x)}{\Omega(x,\alpha(\mu))}.
\]
On account of the identity~\eqref{limit_omega} and the limit~\eqref{limite} we have that 
\[
\lim_{(x,\alpha)\rightarrow(+\infty,-1)}\Omega(x,\alpha)=\lim_{(x,\alpha)\rightarrow(+\infty,-1)}(1+\alpha)\fun(\alpha)\omega(x,\alpha)=+\infty. 
\]The result follows then on account of the continuity of $\alpha$ and due to $\lim_{x\rightarrow+\infty} x\arcsin(M/x)=M$. 
\end{prova}

\begin{lema}\label{lema:omega}
For any $M>0$,
\[
\lim_{(x,\alpha)\rightarrow(+\infty,-1)}\frac{\Omega(x,\alpha)+\mathscr K(\alpha)-\omega(M,\alpha)}{\Omega(x,\alpha)}=1.
\]
\end{lema}
\begin{prova}
The result is straightforward using identity~\eqref{limit_omega} together with the limit~\eqref{limite} and the fact that both $\omega(M,\alpha)$ and $\mathscr K(\alpha)$ are bounded functions for $\alpha\approx -1$.
\end{prova}

The next result requires the introduction of the Gauss hypergeometric series
\[
{}_2F_1(a,b;c;z)=\sum_{n=0}^{\infty}\frac{(a)_n(b)_n}{(c)_n}\frac{z^n}{n!}
\]
for $z,a,b,c\in\C$, where $(a)_n=a(a+1)\cdots(a+n-1)$. The radius of convergence of this series is $1$. In what follows we consider $z,a,b,c\in\R$. The series is absolute convergent when $c-a-b>0$, divergent when $c-a-b\leq -1$ and convergent when $-1<c-a-b\leq 0$ provided that $z\neq 1$. Notice that the series is not defined when $c$ is equal to $-m$, ($m=0,1,2,\dots$), provided $a$ or $b$ is not a negative integer $n$ with $n<m$ (see~\cite{AbrSte92}).

\begin{lema}\label{lema:derivada2F1}
$\frac{d}{dz}{}_2F_1(a,b;a+1;z)=\frac{a}{z}((1-z)^{-b}-{}_2F_1(a,b;a+1;z))$.
\end{lema}

\begin{prova}
This is straightforward by using the formulae in~\cite{AbrSte92}. Indeed, it shows that
\[
\frac{d}{dz}z^a{}_2F_1(a,b;a+1;z)=az^{a-1}{}_2F1(a,b;a;z)=az^{a-1}(1-z)^{-b},
\]
where the first equality is a particular case of $15.2.4$ and the second one follows by applying $15.1.8$. Then an easy manipulation yields to the desired equality after deriving the product on the left.
\end{prova}

The following result is the key tool on the proof of Theorem~\ref{thm:main}.

\begin{prop}\label{prop:cas-facil2}
Let $\Lambda$ be an open subset of $\R^d$, $M>0$ and $\alpha:\Lambda\rightarrow\R$ be a continuous function such that $\alpha(\hat\mu)=-1$ for some fixed $\hat\mu\in\Lambda$. Then 
\[
x^{\alpha(\mu)+1}\int_{\arcsin(M/x)}^{\frac{\pi}{2}}(\sin\theta)^{\alpha(\mu)} d\theta \sim_{\infty} \Omega(x,\alpha(\mu)) \text{ at }\hat\mu.
\]
\end{prop}
\begin{prova}
First we claim that, for all $x>M$,
\[
\int_{\arcsin(M/x)}^{\frac{\pi}{2}}(\sin\theta)^{\alpha} d\theta=
\begin{cases}
\sqrt{1-\frac{M^2}{x^2}}\ {}_2F_1\!\!\left(\frac{1}{2},\frac{1-\alpha}{2};\frac{3}{2};1-\frac{M^2}{x^2}\right) & \text{ if }\alpha\neq -1,\\
\mathrm{arctanh}\!\left(\!\sqrt{1-\frac{M^2}{x^2}}\right) & \text{ if }\alpha=-1.
\end{cases}
\]
In order to prove the claim, let us start considering the case $\alpha=-1$. We can write 
\[
\int_{\arcsin(M/x)}^{\frac{\pi}{2}}\frac{d\theta}{\sin\theta} = \int_{M/x}^{1} \frac{dz}{z\sqrt{1-z^2}}
\]
for any $M>0$ and the result follows using that $\frac{d}{dz}\mathrm{arctanh}(\sqrt{1-z^2})=\frac{1}{z\sqrt{1-z^2}}$. In the case $\alpha\neq -1$, on account of Lemma~\ref{lema:derivada2F1}, we observe that
\[
\frac{d}{d\theta}\left(-{}_2F_1\!\left(\frac{1}{2},\frac{1-\alpha}{2};\frac{3}{2},\cos^2\theta\right)\cos\theta\right)=(\sin\theta)^{\alpha}.
\]
The claim follows, in this second case, evaluating the primitive at the endpoints of the interval of integration using that ${}_2F_1\left(\frac{1}{2},\frac{1-\alpha}{2};\frac{3}{2},0\right)=1$ and $\cos(\arcsin(M/x))=\sqrt{1-M^2/x^2}$. 

Let us now prove the result. To do so, we take advantage of equality~$15.3.6$ on~\cite{AbrSte92}. That is,
\begin{align*}
{}_2F_1(a,b;c;z)&=\frac{\Gamma(c)\Gamma(c-a-b)}{\Gamma(c-a)\Gamma(c-b)}{}_2F_1(a,b;a+b-c+1;1-z)\\&
\phantom{=}+(1-z)^{c-a-b}\frac{\Gamma(c)\Gamma(a+b-c)}{\Gamma(a)\Gamma(b)}{}_2F_1(c-a,c-b;c-a-b+1;1-z)
\end{align*}
for all $z\in\C$ with $\abs{\text{arg}(1-z)}<\pi$ and $c\neq a+b\pm m$, $m\in\N$. In the particular case that we are concerned with, the previously equality yields to
\[
{}_2F_1\!\!\left(\frac{1}{2},\frac{1-\alpha(\mu)}{2};\frac{3}{2};1-\frac{M^2}{x^2}\right)=\fun(\alpha(\mu))\frac{x}{\sqrt{x^2-M^2}}-\frac{M^{1+\alpha(\mu)}}{(1+\alpha(\mu))x^{1+\alpha(\mu)}}{}_2F_1\left(1,1+\frac{\alpha(\mu)}{2};\frac{3+\alpha(\mu)}{2};\frac{M^2}{x^2}\right)
\]
provided that $\alpha(\mu)\neq -1$. Here we used that ${}_2F_1(a,b;b;z)=(1-z)^{-a}$ (see $15.1.8$ in~\cite{AbrSte92}). 
%Thus we have
%\[
%x^{\alpha(\mu)+1}\int_{\arcsin(M/x)}^{\frac{\pi}{2}}(\sin\theta)^{\alpha(\mu)} d\theta=
%\begin{cases}
%x^{\alpha(\mu)+1}\fun(\alpha(\mu))-\frac{M^{1+\alpha(\mu)}}{1+\alpha(\mu)}\sqrt{1-\frac{M^2}{x^2}} {}_2F_1\left(1,1+\frac{\alpha(\mu)}{2};\frac{3+\alpha(\mu)}{2};\frac{M^2}{x^2}\right)& \text{if }\alpha(\mu)\neq -1,\\
%\mathrm{arctanh}\!\left(\!\sqrt{1-\frac{M^2}{x^2}}\right) & \text{if }\alpha(\mu)=-1.
%\end{cases}
%\]

The result will follow once we prove that
\[
x^{\alpha(\mu)+1}\int_{\arcsin(M/x)}^{\frac{\pi}{2}}(\sin\theta)^{\alpha(\mu)}d\theta = \Omega(x,\alpha(\mu))+\mathscr K(\alpha)-\omega(M,\alpha(\mu))+\frac{1}{x^2}(c+r(x,\mu))
\]
for some constant $c$ and some smooth function $r(x,\mu)$ satisfying $\lim_{(x,\mu)\rightarrow (+\infty,\hat\mu)}r(x,\mu)= 0$. Indeed, if it is so, we would have
\[
\frac{x^{\alpha(\mu)+1}}{\Omega(x,\alpha(\mu))}\int_{\arcsin(M/x)}^{\frac{\pi}{2}}(\sin\theta)^{\alpha(\mu)}d\theta = \frac{\Omega(x,\alpha(\mu))+\mathscr K(\alpha)-\omega(M,\alpha(\mu))}{\Omega(x,\alpha(\mu))}+\frac{c+r(x,\mu)}{x^2\Omega(x,\alpha(\mu))}.
\]
Since, by Lemma~\ref{lema:omega},
\[
\lim_{(x,\mu)\rightarrow(+\infty,\hat\mu)}\frac{\Omega(x,\alpha(\mu))+\mathscr K(\alpha(\mu))-\omega(M,\alpha(\mu))}{\Omega(x,\alpha(\mu))}=1,
\]
the result would follow using the limits \eqref{limit_omega} and~\eqref{limite}, which imply $\lim_{(x,\mu)\rightarrow(+\infty,\hat\mu)}x^2\Omega(x,\alpha(\mu))=+\infty$.

We claim at this point that 
\[
\sqrt{1-\frac{M^2}{x^2}}{}_2F_1\left(1,1+\frac{\alpha(\mu)}{2};\frac{3+\alpha(\mu)}{2};\frac{M^2}{x^2}\right)=1+\frac{1+\alpha(\mu)}{x^2}\left(\frac{M^2}{2(3+\alpha(\mu))}+r(x,\mu)\right),
\]
with $\lim_{(x,\mu)\rightarrow(+\infty,\hat\mu)}r(x,\mu)=0$. Indeed, from equality $15.3.3$ in~\cite{AbrSte92} we have that
\[
\sqrt{1-\frac{M^2}{x^2}}{}_2F_1\left(1,1+\frac{\alpha(\mu)}{2};\frac{3+\alpha(\mu)}{2};\frac{M^2}{x^2}\right)={}_2F_1\left(\frac{1}{2},\frac{1+\alpha(\mu)}{2};\frac{3+\alpha(\mu)}{2};\frac{M^2}{x^2}\right).
\]
By definition,
\[
{}_2F_1\left(\frac{1}{2},\frac{1+\alpha(\mu)}{2};\frac{3+\alpha(\mu)}{2};\frac{M^2}{x^2}\right)=1+\frac{(1+\alpha(\mu))}{2(3+\alpha(\mu))}\frac{M^2}{x^2}+\sum_{n=2}^{\infty}\frac{\bigl(\frac{1}{2}\bigr)_{\!n}\bigl(\frac{1+\alpha(\mu)}{2}\bigr)_{\!n} M^{2n}}{\bigl(\frac{3+\alpha(\mu)}{2}\bigr)_{\!n} n!}x^{-2n}.
\]
Notice that, in the previous expression, there is always a factor $(1+\alpha(\mu))$ on the numerator for $n\geq 1$. Moreover, it is easy to show that, after extracting $(1+\alpha(\mu))$ as a common factor, the sequence has positive monotone decreasing coefficients for $\alpha(\mu)\approx -1$. This follows from the identity $(a)_{n+1}=(a_n)(a+n)$. Therefore,
\[
\frac{1}{2(3+\alpha(\mu))}\frac{M^2}{x^2}\leq \frac{{}_2F_1\left(\frac{1}{2},\frac{1+\alpha(\mu)}{2};\frac{3+\alpha(\mu)}{2};\frac{M^2}{x^2}\right)-1}{1+\alpha(\mu)}\leq\frac{1}{2(3+\alpha(\mu))}\left(\frac{M^2}{x^2}+\sum_{n=2}^{\infty}M^{2n}x^{-2n}\right).
\]
The claim follows then on account that the series $\sum_{n=2}^{\infty}z^{2n}$ is convergent for $\abs{z}<1$ and taking the limit as $(x,\mu)\rightarrow(+\infty,\hat\mu)$. Here we point out that the constant $M$ is fixed.

On account of the previous claim, after some elementary manipulations, we have, for $\alpha(\mu)\neq -1$,
\begin{align*}
x^{1+\alpha(\mu)}\int_{\arcsin(M/x)}^{\frac{\pi}{2}}(\sin\theta)^{\alpha(\mu)}d\theta&=
x^{\alpha(\mu)+1}\fun(\alpha(\mu))-\frac{M^{1+\alpha(\mu)}}{1+\alpha(\mu)}-\frac{1}{x^2}\left(\frac{M^{3+\alpha(\mu)}}{2(3+\alpha(\mu))}+r(x,\mu)\right)\\ 
&=
\Omega(x,\alpha(\mu))+\mathscr K(\alpha(\mu))-\omega(M,\alpha(\mu))-\frac{1}{x^2}\left(\frac{M^{3+\alpha(\mu)}}{2(3+\alpha(\mu))}+r(x,\mu)\right).
\end{align*}
In the case $\alpha(\mu)=-1$, by the Taylor's series expansion,
\begin{align*}
\int_{\arcsin(M/x)}^{\frac{\pi}{2}}(\sin\theta)^{\alpha(\mu)}d\theta&=\mathrm{arctanh}\!\left(\!\!\sqrt{1-\frac{M^2}{x^2}}\right)\\
&=\log(2x)-\log(M)-\frac{1}{x^2}\left(\frac{M^2}{4}+r(x)\right)\\
&=\Omega(x,\alpha(\mu))+\mathscr K(\alpha(\mu))-\omega(M,\alpha(\mu))-\frac{1}{x^2}\left(\frac{M^2}{4}+r(x)\right)
\end{align*}
with $\lim_{x\rightarrow+\infty}r(x)= 0$. Therefore, for any $\alpha(\mu)$, we have
\[
x^{1+\alpha(\mu)}\int_{\arcsin(M/x)}^{\frac{\pi}{2}}(\sin\theta)^{\alpha(\mu)}d\theta=
\Omega(x,\alpha(\mu))+\mathscr K(\alpha(\mu))-\omega(M,\alpha(\mu))-\frac{1}{x^2}\left(\frac{M^{3+\alpha(\mu)}}{2(3+\alpha(\mu))}+r(x,\mu)\right)
\]
with $\lim_{(x,\mu)\rightarrow(+\infty,\hat\mu)}r(x,\mu)=0$. The result follows then by Lemma~\ref{lema:omega}.
\end{prova}

\begin{defi}
Let $\Lambda$ be an open subset of $\R^d$ and $\{f_{\mu}\}_{\mu\in\Lambda}$ be a continuous family of analytic functions on $[0,+\infty)$. Setting $f_0(x,\mu)\!:=f_{\mu}(x)$, we define
$
f_m(x,\mu)\!:=f_{m-1}(x,\mu)x^2+\int_0^x f_{m-1}(s,\mu)ds \text{ for all }m\geq 1.$
\end{defi}

\begin{lema}[see~\cite{ManRojVil2015}]\label{lema:puja}
Let $\Lambda$ be an open subset of $\R^d$ and $\{f_{\mu}\}_{\mu\in\Lambda}$ be a continuous family of analytic functions on $[0,+\infty)$. Then, for any $m\in\N$,
\[
\F[f_{\mu}](x)=\frac{1}{x^{2m}}\F[f_m(\cdot,\mu)](x) \text{ for all }x>0.
\]
\end{lema}

Next result was proved in~\cite[Proposition 2.16]{ManRojVil2015}. We point out that in its original statement the requirement $\alpha(\hat\mu)\neq -2m-1$ is not needed. This requirement was assumed for $b_m(\alpha(\mu))=\prod_{i=1}^m \frac{\alpha(\mu)+2i}{\alpha(\mu)+2i-1}$ to be well defined in a neighborhood of $\hat\mu$. However, the reader may easily check that this is true for $\alpha(\hat\mu)=-2m-1$.

\begin{prop}\label{prop:antiga}
Let $\Lambda$ be an open subset of $\R^d$ and consider a continuous family $\{f_{\mu}\}_{\mu\in\Lambda}$ of analytic functions on $[0,+\infty)$ such that $f_{\mu}(x)\sim_{\infty} a(\mu)x^{\alpha(\mu)}$ at any $\mu\in\Lambda$.  Assume that for some $\hat\mu\in\Lambda$, $\alpha(\hat\mu)\leq -1$ and let $m\in\N\cup\{0\}$ be such that $\alpha(\hat\mu)\in[-2m-1,-2m+1).$ If 
$M_1[f_{\mu}]\equiv M_2[f_{\mu}]\equiv\cdots\equiv M_m[f_{\mu}]\equiv 0$ for all $\mu\in\Lambda$ and $\alpha(\hat\mu)+2m\neq 0$ then $f_m(x,\mu)\sim_{\infty} a(\mu)b_m(\alpha(\mu))x^{\alpha(\mu)+2m}$.
\end{prop}

\begin{prooftext}{Proof of Theorem~\ref{thm:main}.}
Let us fix $m\in\N\cup\{0\}$ such that $\alpha(\hat\mu)+2m=-1$ and let us assume that $M_1[f_{\mu}]\equiv M_2[f_{\mu}]\equiv \cdots M_m[f_{\mu}]\equiv 0$. Moreover, let us set $\hat a\!:=a(\hat\mu)b_m(\alpha(\hat\mu))$ for the sake of shortness. The result will follow once we show that for any $\epsilon>0$ there exist $x_0>0$ and $\delta_0>0$ such that 
\[
\abs{\frac{x}{\Omega(x,\alpha(\mu)+2m)}\int_0^{\frac{\pi}{2}}f_{m}(x\sin\theta,\mu)d\theta-\hat a}<\epsilon
\]
for all $x>x_0$ and $\|\mu-\hat\mu\|<\delta_0$. Indeed, if it is so we have that 
\[
x\F[f_{m}(\cdot,\mu)](x)\sim_{\infty} a(\mu)b_m(\alpha(\mu))\Omega(x,\alpha(\mu)+2m) \text{ at }\hat\mu.
\]
Consequently, on account of Lemma~\ref{lema:puja}, $x^{2m+1}\F[f_{\mu}](x)\sim_{\infty} a(\mu)b_m(\alpha(\mu))\Omega(x,\alpha(\mu)+2m) \text{ at }\hat\mu,$ as we desired.

By Proposition~\ref{prop:antiga} we have that $f_m(x,\mu)\sim_{\infty} a(\mu)b_m(\alpha(\mu))x^{\alpha(\mu)+2m}$ in a neighborhood of $\hat\mu$ and $\alpha(\hat\mu)+2m=-1$. Then, for any fixed $\epsilon>0$, there exist $M>0$ and $\delta_1>0$ such that 
\begin{equation}\label{f_eq}
\abs{x^{-\alpha(\mu)-2m}f_{m}(x,\mu)-\hat a}<\frac{\epsilon}{8}
\end{equation}
for all $x>M$ and $\|\mu-\hat\mu\|<\delta_1$. Moreover we can assume that $\alpha(\mu)+2m\in(-2,0)$ for all $\|\mu-\hat\mu\|<\delta_1$.
Thus, by Lemma~\ref{lema:previ}, there exists $x_1>M$ such that
\begin{equation}\label{f_ineq1}
\abs{\frac{x}{\Omega(x,\alpha(\mu)+2m)}\int_0^{\arcsin(M/x)}f_m(x\sin\theta,\mu)d\theta}\leq \frac{\epsilon}{2}
\end{equation}
for all $x>x_1$.
From Proposition~\ref{prop:cas-facil2} we have that there exist $x_0>x_1$ and $0<\delta_0<\delta_1$ such that
\begin{equation}\label{f_ineq2}
\abs{\frac{\hat a x}{\Omega(x,\alpha(\mu)+2m)}\int_{\arcsin(M/x)}^{\frac{\pi}{2}}(x\sin\theta)^{\alpha(\mu)+2m}d\theta-\hat a}\leq \frac{\epsilon}{4}
\end{equation}
and
\begin{equation}\label{f_ineq3}
\abs{\frac{x}{\Omega(x,\alpha(\mu)+2m)}\int_{\arcsin(M/x)}^{\frac{\pi}{2}}(x\sin\theta)^{\alpha(\mu)+2m}d\theta}\leq 2
\end{equation}
for all $x>x_0$ and $\|\mu-\hat\mu\|<\delta_0$.
Lastly, for all $x>x_0$ and $\|\mu-\hat\mu\|<\delta_0$,
\begin{align*}
\abs{\frac{x}{\Omega(x,\alpha(\mu)+2m)}\int_{\arcsin(M/x)}^{\frac{\pi}{2}}\left(\frac{f_{m}(x\sin\theta,\mu)}{(x\sin\theta)^{\alpha(\mu)+2m}}-\hat a\right)(x\sin\theta)^{\alpha(\mu)+2m}d\theta}&\leq\\
&\hspace*{-5cm}\leq \frac{\epsilon}{8} \frac{x}{\Omega(x,\alpha(\mu)+2m)}\int_{\arcsin(M/x)}^{\frac{\pi}{2}}(x\sin\theta)^{\alpha(\mu)+2m}d\theta\leq \frac{\epsilon}{2},
\end{align*}
where we used~\eqref{f_eq} in the first inequality and~\eqref{f_ineq3} in the second one. The previous inequality together with~\eqref{f_ineq1} and~\eqref{f_ineq2} proves the Theorem.
\end{prooftext}

\subsection{Proof of Theorem~\ref{thm:main-gen}}

In this section we prove Theorem~\ref{thm:main-gen} as a corollary of Theorem~\ref{thm:main}. We first recall a result that was already shown in~\cite{ManRojVil2016-JDDE}.

\begin{prop}\label{com}
Let $f\in C^{\omega}(0,+\infty)$. If $f$ can be extended analytically to $x=0$, then $\LL_{\boldsymbol\nu_{n}}[f]$ can be extended analytically to $x=0$. Moreover, $\F\circ\LL_{\boldsymbol\nu_n}=\LL_{\boldsymbol\nu_n}\circ \F.$
\end{prop}

\begin{prooftext}{Proof of Theorem~\ref{thm:main-gen}.}
By Proposition~\ref{com} $\LL_{\boldsymbol\nu_{n}(\mu)}[f_{\mu}]$ is an analytic function on $[0,+\infty)$ for each $\mu\in\Lambda$ and
\[
(\LL_{\boldsymbol\nu_{n}(\mu)}\circ\F)[f_{\mu}](x)=(\F\circ \LL_{\boldsymbol\nu_{n}(\mu)})[f_{\mu}](x)=\int_0^{\frac{\pi}{2}}\LL_{\boldsymbol\nu_{n}(\mu)}[f_{\mu}](x\sin\theta)d\theta.
\]
Then the result follows by applying Theorem~\ref{thm:main} to the family $\{\LL_{\boldsymbol\nu_{n}(\mu)}[f_{\mu}]\}_{\mu\in\Lambda}$.
\end{prooftext}

\section{Criticality of the period function at the outer boundary}\label{sec:criticality}

This section is devoted to the dynamical results of the paper. Consider the family of analytic potential differential systems~\eqref{sys} depending on a parameter $\mu\in\Lambda\subset\R^d$ and suppose that the origin is a non-degenerate center for all $\mu$. We denote the projection of the period annulus on the $x$-axis by $\PI_{\mu}=(x_{\ell}(\mu),x_r(\mu))$, $x_{\ell}(\mu)<0<x_r(\mu),$ and  by $h_0(\mu)$ the energy level at the outer boundary of the period annulus.

\begin{defi}\label{def:Ak}
We say that the family of systems~\eqref{sys} verifies the hypothesis \hypertarget{H}{\hyperlink{H}{\textnormal{\textbf{(H)}}}} in case that:
\begin{enumerate}[$(a)$]
\item For all $k\geq 0$, the map $(x,\mu)\longmapsto V_{\mu}^{(k)}(x)$ is continuous on $\{(x,\mu)\in\R\times\Lambda : x\in I_{\mu}\},$
\item $\mu\longmapsto x_{r}(\mu)$ is continuous on $\Lambda$ or $x_{r}(\mu)=+\infty$ for all $\mu\in\Lambda,$
\item $\mu\longmapsto x_{\ell}(\mu)$ is continuous on $\Lambda$ or $x_{\ell}(\mu)=-\infty$ for all $\mu\in\Lambda,$
\item $\mu\longmapsto h_0(\mu)$ is continuous on $\Lambda$ or $h_0(\mu)=+\infty$ for all $\mu\in\Lambda.$
\end{enumerate}\vspace*{-.5cm}
\end{defi}

Next result is proved in \cite{ManRojVil2015}. We recall that $g_{\mu}(x)=\text{sgn}(x)\sqrt{V_{\mu}(x)}$.

\begin{lema}\label{lema:g_menys_1_continua}
Assume that the family of systems~\eqref{sys} verifies hypothesis \hyperlink{H}{\textnormal{\textbf{(H)}}}.  
Then the map $(x,\mu)\longmapsto g_{\mu}^{-1}(x)$ is continuous on the open set
$\{(x,\mu)\in\R\times\Lambda: x\in
(-\sqrt{h_{0}(\mu)},\sqrt{h_{0}(\mu)})\}$.
\end{lema}

The following sections are concerned with sufficient conditions to bound the criticality at the outer boundary for potential systems~\eqref{sys} verifying hypothesis \hyperlink{H}{\textnormal{\textbf{(H)}}}. \secc{section:infinite} deals with the case $h_0\equiv +\infty$ (see Theorem~\ref{thm:criticality_infinite}), whereas \secc{section:finite} tackles the case $h_0$ finite (see Theorem~\ref{thm:criticality_finite}).

\subsection{Potential systems with infinite energy}\label{section:infinite}

In this section let us consider that the energy at the outer boundary is $h_0(\mu)=+\infty$ for all $\mu\in\Lambda$. Following the strategy in~\cite{ManRojVil2016-JDDE}, we plan to find sufficient conditions such that $fT_{\mu}'$ can be embedded into the ECT-system $(h^{\nu_1(\mu)},h^{\nu_2(\mu)},\dots,h^{\nu_n(\mu)})$, where $f$ is an analytic non-vanishing function. Next result is analogous to \cite[Lemma 3.5]{ManRojVil2016-JDDE} where only the power $h^{\nu_n(\mu)}$ appears multiplying the Wronskian. The reader may check that the same proof there can be applied here, so we skip it for the sake of shortness.

\begin{lema}\label{lema:crit_Wronskia_infinit}
Let $\{X_{\mu}\}_{\mu\in\Lambda}$ be a family of potential analytic differential systems verifying~\hyperlink{H}{\textnormal{\textbf{(H)}}} and such that $h_0\equiv +\infty$. Assume that there exist $n\geqslant 1$ continuous functions $\nu_1,\nu_2\dots,\nu_n$ in a neighborhood of some fixed $\hat\mu\in\Lambda$, a continuous function $\alpha:\Lambda\rightarrow\R$ with $\alpha(\hat\mu)=-1$ and an analytic non-vanishing function $f$ on $(0,+\infty)$ such that
\[
\lim_{(h,\mu)\rightarrow(+\infty,\hat\mu)}\frac{h^{\nu_n(\mu)}}{\Omega(h,\alpha(\mu))}W[h^{\nu_1(\mu)},\dots,h^{\nu_{n-1}(\mu)},f(h)T_{\mu}'(h)]=\ell\neq 0.
\]
Then $\mathrm{Crit}\bigl((\Pi_{\hat\mu},X_{\hat\mu}),X_{\mu}\bigr)\leq n-1$.
\end{lema}

We are now in position to state the main result concerning the criticality at the outer boundary for the case $h_0\equiv +\infty.$ In its statement, and from now on, for a given function $f\!:\!(-a,a)\longrightarrow\R$, $a>0$, we denote $\mathcal P[f](x)\!:=f(x)+f(-x)$. Let us also remark that the assumption requiring the existence of functions $\nu_1,\nu_2,\ldots,\nu_n$ is void in case that $n=0.$ The same happens for the assumption of $M_1\equiv M_2\equiv \cdots\equiv M_m\equiv 0$ when $m=0$.

\begin{thmx}\label{thm:criticality_infinite}
Let $\{X_{\mu}\}_{\mu\in\Lambda}$ be a family of potential analytic differential systems verifying~\hyperlink{H}{\textnormal{\textbf{(H)}}} with $h_0\equiv+\infty$ and that there exist $n\geq 0$ continuous functions $\nu_1,\nu_2,\ldots,\nu_n$ in a neighborhood of some fixed $\hat\mu\in\Lambda$ such that 
\[
(\LL_{\boldsymbol\nu_n(\mu)}\circ\mathcal P)[z(g_{\mu}^{-1})''(z)]\sim_{\infty} a(\mu)x^{\xi(\mu)} \text{ at }\hat\mu.
\]
For each $i\in\N$, let $M_i(\mu)$ be the $i$-th momentum of $(\LL_{\boldsymbol\nu_n(\mu)}\circ\mathcal P)[z(g_{\mu}^{-1})''(z)]$, whenever it is well defined. If $\xi(\hat\mu)=-1-2m$ for some $m\in\N\cup\{0\}$ and $M_1\equiv M_2\equiv\cdots\equiv M_m\equiv 0$ then 
$\mathrm{Crit}\bigl((\out_{\hat\mu},X_{\hat\mu}),X_{\mu}\bigr)\leqslant n.$
\end{thmx}

\begin{prova}
Let us denote $f_{\mu}(z)\!:=\PP[z\bigl(g^{-1}_{\mu}\bigr)''(z)]$. By \lemc{lema:g_menys_1_continua} and the hypothesis \hyperlink{H}{\textnormal{\textbf{(H)}}} we have that $\{f_{\mu}\}_{\mu\in\Lambda}$ is a continuous family of analytic functions on $(0,+\infty)$ that extends analytically to $z=0$. Since $\xi(\mu)$ is the quantifier of $\{\LL_{\boldsymbol\nu_n(\mu)}[f_{\mu}]\}_{\mu\in\Lambda}$ at infinity, $M_1\equiv M_2\equiv\cdots \equiv M_m\equiv 0$ and $\xi(\hat\mu)+2m=-1$, applying Theorem~\ref{thm:main-gen} we can assert that 
\[
x^{2m+1}(\LL_{\boldsymbol\nu_n(\mu)}\circ\F)[f_{\mu}](x)\sim_{\infty} a(\mu)b_m(\mu)\Omega(x,\xi(\mu)+2m) \text{ at }\hat\mu.
\] 
Then, on account of the definition of $\LL_{\boldsymbol\nu_n(\mu)}$, see~\refc{eq:op_Wronskia}, we have that
\[
\lim_{(h,\mu)\rightarrow(+\infty,\hat\mu)}\frac{h}{\Omega(h,\xi(\mu)+2m)}\frac{W\bigl[h^{\nu_1(\mu)},\dots,h^{\nu_{n}(\mu)},\F[f_{\mu}](h)\bigr]}{h^{\sum_{i=1}^{n}(\nu_i(\mu)-i)}}\neq 0.
\]
The result follows then by \lemc{lema:crit_Wronskia_infinit} and using equality~\eqref{formula}.
\end{prova}

\subsection{Potential systems with finite energy}\label{section:finite}

In this section let us consider that the energy at the outer boundary is finite for all $\mu\in\Lambda$. With the intention of embedding $fT_{\mu}'$ into some ECT-system for an appropriate non-vanishing function $f$, we proceed as in~\cite{ManRojVil2016-JDDE} and ``translate'' the case $h_0<+\infty$ to the case $h_0=+\infty$ so we can take advantage of Theorem~\ref{thm:main-gen}. With this aim in view, we define next a differential operator which is conjugated to $\LL_{\boldsymbol\nu_n(\mu)}$. The conjugation is precisely the tool that enables this translation. Given $\nu_1,\dots,\nu_n\in\R$, consider the linear ordinary differential operator 
\[
\DD_{\boldsymbol\nu_n}\!: C^{\omega}(0,1)\longrightarrow C^{\omega}(0,1)
\]
defined by
\begin{equation}\label{eq:D}
\DD_{\boldsymbol\nu_n}[f](x)\!:=(x(1-x^2))^{\frac{n(n+1)}{2}}\frac{W\left[\psi_{\nu_1},\ldots,\psi_{\nu_n},f\right](x)}{\prod_{i=1}^n\psi_{\nu_i}(x)},
\end{equation}
where we use the notation $\boldsymbol\nu_n=(\nu_1,\dots,\nu_n)$ and 
$ \psi_{\nu}(x)\!:=\frac{1}{1-x^2}\left(\frac{x}{\sqrt{1-x^2}}\right)^{\nu}.$
Furthermore we define $\DD_{\boldsymbol\nu_0}\!:=id$ for the sake of completeness. Setting
$\phi(x)\!:=\frac{x}{\sqrt{1+x^2}},$ 
we also consider the operator 
\[
\CC\!:C^{\omega}[0,1)\longrightarrow C^{\omega}[0,+\infty)\]
defined by
\begin{equation}\label{defi:operator_C}
\CC[f](x)\!:=\bigl(1-\phi^2(x)\bigr)\bigl(f\circ \phi\bigr)(x)=\frac{1}{1+x^2}\bigl(f\circ \phi\bigr)(x).
\end{equation}
This operator is the conjugation mentioned before.

\begin{defi}
Let $f\in C^{\omega}[0,1).$ Then, for each $n\in\N$, we call
\[
N_n[f]:=\int_0^{1} \frac{f(x)}{\sqrt{1-x^2}}\left(\frac{x}{\sqrt{1-x^2}}\right)^{2n-2}dx
\]
the $n$-th \emph{momentum} of $f$, whenever it is well defined.
\end{defi}

Next results show the way $\CC$ conjugates $\DD_{\boldsymbol\nu_n}$ and $\LL_{\boldsymbol\nu_n}$. We refer the reader to~\cite{ManRojVil2016-JDDE} for the proofs.

\begin{lema}\label{lem:lemes_varis}
Consider $\nu_1,\nu_2,\dots,\nu_n\in\R$. Then the following hold:
\begin{enumerate}[$(a)$]
\item $\CC[\psi_{\nu_i}](x)=x^{\nu_i}$ for $i=1,2,\ldots,n.$
\item $\CC\circ\DD_{\boldsymbol\nu_n}=\LL_{\boldsymbol\nu_n}\!\circ\CC$.
\item $\bigl(\F\circ \CC\bigr)[f](x)=\sqrt{1+x^2}\bigl(\CC\circ \F\bigr)[f](x)$ for any $f\in C^{\omega}(0,1).$ 
\item $N_n=M_n\circ \CC$.
\end{enumerate}
\end{lema}

\begin{lema}\label{lem:traduccio_quantificador}
Let $\{f_{\mu}\}_{\mu\in\Lambda}$ be a continuous family of analytic
functions on~$[0,1)$. Then 
\[
f_{\mu}(z)\sim_{z=1}a(\mu)(1-z)^{-\alpha(\mu)}\text{ at }\hat\mu\text{ if and only if }\CC[f_{\mu}](x)\sim_{\infty} a(\mu)2^{\alpha(\mu)}x^{2\alpha(\mu)-2}\text{ at }\hat\mu.
\]
\end{lema}

The next result is well known (see \cite{ManVil11}).

\begin{lema}\label{lem:Wronskians}
Let $f_0,f_1,\dots,f_{n-1}$ be analytic functions. Then the following statements hold:
\begin{enumerate}[$(a)$]
\item $W[f_0\circ \varphi,\dots, f_{n-1}\circ \varphi](x)=\left(\varphi'(x)\right)^{\frac{(n-1)n}{2}}W[f_0,\dots,f_{n-1}](\varphi(x))$ for any analytic diffeomorphism $\varphi$.
\item $W[gf_0,\dots, gf_{n-1}](x)=g(x)^n W[f_0,\dots,f_{n-1}](x)$ for any analytic function $g$.
\end{enumerate}
\end{lema}

The following Lemma is analogous to Lemma~\ref{lema:crit_Wronskia_infinit} for the case $h_0$ finite. As before, the reader may check that the proof in~\cite[Lemma 3.10]{ManRojVil2016-JDDE} is also valid for this generalization.

\begin{lema}\label{lema:crit_Wronskia_finit}
Let $\{X_{\mu}\}_{\mu\in\Lambda}$ be a family of potential analytic differential systems verifying~\hyperlink{H}{\textnormal{\textbf{(H)}}} and such that $\mu\longmapsto h_0(\mu)$ is continuous on $\Lambda$. Assume that there exist $n\geqslant 1$ continuous functions $\nu_1,\nu_2\dots,\nu_n$ in a neighborhood of some fixed $\hat\mu\in\Lambda$, a continuous function $\alpha:\Lambda\rightarrow\R$ with $\alpha(\hat\mu)=-1$ and an analytic non-vanishing function $f$ on $(0,1)$ such that
\[
\lim_{(z,\mu)\rightarrow (1,\hat\mu)}\frac{(1-z)^{\nu_n(\mu)}}{\Omega\left(\frac{z}{\sqrt{1-z^2}},\alpha(\mu)\right)}W\bigl[\psi_{\nu_1(\mu)}(z),\dots,\psi_{\nu_{n-1}(\mu)}(z),f(z)T_{\mu}'(z^2h_0(\mu))\bigr]=\ell\neq 0.
\]
Then $\mathrm{Crit}\bigl((\Pi_{\hat\mu},X_{\hat\mu}),X_{\mu}\bigr)\leq n-1$.
\end{lema}

Next we state the main result to bound the criticality at the outer boundary of the family of systems~\eqref{sys} in case that its energy level is finite. Again we stress that the assumptions requiring the existence of functions $\nu_1,\nu_2,\ldots,\nu_n$ for $n=0$ and $N_1\equiv N_2\equiv\cdots\equiv N_m\equiv 0$ for $m=0$ are void in its statement.

\begin{thmx}\label{thm:criticality_finite}
Let us assume that the family of potential analytic systems~\eqref{sys} verifies hypothesis \hyperlink{H}{\textnormal{\textbf{(H)}}} with $h_0(\mu)<+\infty$ for all $\mu\in\Lambda$ and that there exist $n\geq 0$ continuous functions $\nu_1,\nu_2,\ldots,\nu_n$ in a neighborhood of some fixed $\hat\mu\in\Lambda$ such that 
\[
(\DD_{\boldsymbol\nu_n(\mu)}\circ\mathcal P)\bigl[z\sqrt{h_0(\mu)}(g_{\mu}^{-1})''(z\sqrt{h_0(\mu)})\bigr]\sim_{z=1} a(\mu)(1-z)^{\xi(\mu)}\text{ at }\hat\mu.
\]
For each $i\in\N$, let $N_i(\mu)$ be the $i$-th momentum of $(\DD_{\boldsymbol\nu_n(\mu)}\circ\mathcal P)\bigl[z\sqrt{h_0(\mu)}(g_{\mu}^{-1})''(z\sqrt{h_0(\mu)})\bigr]$, whenever it is well defined. If $\xi(\hat\mu)=\frac{1}{2}-m$ for some $m\in\N\cup\{0\}$ and $N_1\equiv N_2\equiv \ldots\equiv N_m\equiv 0$ then $\mathrm{Crit}\bigl((\out_{\hat\mu},X_{\hat\mu}),X_{\mu}\bigr)\leqslant~n$.
\end{thmx}

\begin{prova}
Let us denote $f_{\mu}(z)\!:=\mathcal P[z\sqrt{h_0(\mu)}(g_{\mu}^{-1})''(z\sqrt{h_0(\mu)})]$. By Lemma~\ref{lema:g_menys_1_continua} and the hypothesis \hyperlink{H}{\textnormal{\textbf{(H)}}} we have that $\{f_{\mu}\}_{\mu\in\Lambda}$ is a continuous family of analytic functions on $(0,1)$ that extends analytically to $z=1$. 
The identity~\eqref{formula}, after the appropriate rescaling, yields to the identity
\begin{equation}\label{eq:dT_operator_finit}
\F[f_{\mu}](z)=\sqrt{2}h_0(\mu)z^2T_{\mu}'(h_0(\mu)z^2) \text{ for all $z\in (0,1).$}
\end{equation}
On account of this equality, to prove the result we must show that there exist $\varepsilon>0$ and a neighborhood $U$ of $\hat\mu$ such that $\F[f_{\mu}](z)$ has at most~$n$ zeros for $z\in (1-\varepsilon,1)$, multiplicities taking into account, for all $\mu\in U.$ By hypothesis $\xi(\hat\mu)=\frac{1}{2}-m$, for some $m\in\N\cup\{0\}$, is the quantifier of $(\DD_{\boldsymbol\nu_n(\mu)}\circ\mathcal P)\bigl[z\sqrt{h_0(\mu)}(g_{\mu}^{-1})''(z\sqrt{h_0(\mu)})\bigr]$ at $z=1$ in $\mu=\hat\mu$. By $(b)$ and $(d)$ in \lemc{lem:lemes_varis},
\[
 N_{i}\bigl[\DD_{\boldsymbol\nu_n(\mu)}[f_{\mu}]\bigr]=
 M_{i}\bigl[(\CC\circ\DD_{\boldsymbol\nu_n(\mu)})[f_{\mu}]\bigr]=
 M_{i}\bigl[(\LL_{\boldsymbol\nu_n(\mu)}\circ\CC)[f_{\mu}]\bigr]
\]
for all $i=1,2,\dots,m$. Therefore the condition $N_1\equiv N_2\equiv\cdots N_m\equiv 0$ in the statement is equivalent to the condition $M_1\equiv M_2\equiv\cdots M_m\equiv 0$, where $M_i$ is the $i$-th momentum of $(\LL_{\boldsymbol\nu_n(\mu)}\circ\CC)[f_{\mu}]$.
Recall at this point that, by $(b)$ in \lemc{lem:lemes_varis}, $\CC\circ\DD_{\boldsymbol\nu_n}=\LL_{\boldsymbol\nu_n}\circ\CC$. By assumption, $\DD_{\boldsymbol\nu_n(\mu)}[f_{\mu}](z)\sim_{z=1}a(\mu)(1-z)^{\xi(\mu)}$ at $\hat\mu$. Then applying~\lemc{lem:traduccio_quantificador} we have that 
\[
(\LL_{\boldsymbol\nu_n(\mu)}\circ\CC)[f_{\mu}](x)\sim_{\infty} a(\mu)2^{\xi(\mu)}x^{\eta(\mu)}\text{ at }\hat\mu,
\] 
where $\eta(\mu)\!:=2\xi(\mu)-2$. Notice that $\eta$ satisfies $\eta(\hat\mu)+2m=-1$. Therefore Theorem~\ref{thm:main-gen} applied to the family $\{(\LL_{\boldsymbol\nu_n(\mu)}\circ\CC)[f_{\mu}]\}_{\mu\in\Lambda}$ shows that
\[
x^{2m+1}(\LL_{\boldsymbol\nu_n(\mu)}\circ\F\circ\CC)[f_{\mu}](x)\sim_{\infty}a(\mu)2^{\xi(\mu)}b_m(\eta(\mu))\Omega(x,\eta(\mu)+2m)\text{ at }\hat\mu
\] 
with $a(\hat\mu)b_m(\eta(\hat\mu))\neq 0$.
Let us note that 
\begin{align*}
(\LL_{\boldsymbol\nu_n(\mu)}\circ\F\circ\CC)[f_{\mu}](x)&=\LL_{\boldsymbol\nu_n(\mu)}\bigl[\sqrt{1+x^2}(\CC\circ\F)[f_{\mu}](x)\bigr]\\
&=(\LL_{\boldsymbol\nu_n(\mu)}\circ\CC)\bigl[(1-z^2)^{-\frac{1}{2}}\F[f_{\mu}](z)\bigr](x)\\
&=(\CC\circ\DD_{\boldsymbol\nu_n(\mu)})\bigl[(1-z^2)^{-\frac{1}{2}}\F[f_{\mu}](z)\bigr](x),
\end{align*}
with $z=\phi(x)=\frac{x}{\sqrt{1+x^2}}$, where we use $(c)$ in \lemc{lem:lemes_varis} in the first equality, the identity $\sqrt{1+x^2}\CC[\varphi](x)=\CC[(1-z^2)^{-\frac{1}{2}}\varphi(z)]$
with $\varphi=\F[f_{\mu}]$ in the second one, and $(b)$ in \lemc{lem:lemes_varis} in the third one. Note also that $\bigl\{(\LL_{\boldsymbol\nu_n(\mu)}\circ\F\circ\CC)[f_{\mu}]\bigr\}_{\mu\in\Lambda}$ is a continuous family of analytic functions on $[0,+\infty).$ On account of the definition of $\CC$ in~\eqref{defi:operator_C} and the previous equality, we have that
\[
\frac{x^{2m+1}}{1+x^2}\DD_{\boldsymbol\nu_n(\mu)}\bigl[(1-\phi(x)^2)^{-\frac{1}{2}}\F[f_{\mu}](\phi(x))\bigr]\sim_{\infty}a(\mu)2^{\xi(\mu)}b_m(\eta(\mu))\Omega(x,\eta(\mu)+2m) \text{ at }\hat\mu.
\]
Setting $x=\phi^{-1}(z)=\frac{z}{\sqrt{1-z^2}}$, the previous identity yields to
\[
\lim_{(z,\mu)\rightarrow(1,\hat\mu)}
\frac{(1-z)^{-2m}}{\Omega\left(\frac{z}{\sqrt{1-z^2}},\eta(\mu)+2m\right)}\DD_{\boldsymbol\nu_n(\mu)}\bigl[(1-z^2)^{-\frac{1}{2}}\F[f_{\mu}](z)\bigr]\neq 0.
\]
Thus, on account of the definition of $\DD_{\boldsymbol\nu_n(\mu)}$ in~\eqref{eq:D},
\[
\lim_{(z,\mu)\rightarrow(1,\hat\mu)}
\frac{(1-z)^{-2m}}{\Omega\left(\frac{z}{\sqrt{1-z^2}},\eta(\mu)+2m\right)}(z(1-z^2))^{\frac{n(n+1)}{2}}\frac{W\bigl[\psi_{\nu_1(\mu)}(z),\dots,\psi_{\nu_n(\mu)}(z),(1-z^2)^{-\frac{1}{2}}\F[f_{\mu}](z)\bigr]}{\prod_{i=1}^n \psi_{\nu_i(\mu)}(z)}\neq 0,
\]
which, since $\psi_{\nu}(z)=\frac{1}{1-z^2}\left(\frac{z}{\sqrt{1-z^2}}\right)^{\nu}$, implies that
\[
\lim_{(z,\mu)\rightarrow(1,\hat\mu)}
\frac{(1-z)^{\kappa(\mu)}}{\Omega\left(\frac{z}{\sqrt{1-z^2}},\eta(\mu)+2m\right)}W\bigl[\psi_{\nu_1(\mu)}(z),\dots,\psi_{\nu_n(\mu)}(z),(1-z^2)^{-\frac{1}{2}}\F[f_{\mu}](z)\bigr]\neq 0
\]
with $\kappa(\mu)\!:=-2m+\frac{n(n+3)}{2}+\frac{1}{2}\sum_{i=1}^n\nu_i(\mu)$. The result follows then by \lemc{lema:crit_Wronskia_finit} and taking the identity \refc{eq:dT_operator_finit} into account. 
\end{prova}

\section{Applications}

\subsection{Proof of Theorem~\ref{thm:familia}}

In this section we resume the study that we began in~\cite{ManRojVil2015,ManRojVil2016-JDDE} for the family of potential differential systems~\eqref{Xmu2} with $\mu=(q,p)\in\Lambda\!:=\{(q,p)\in\R^2:p>q\}$. Following the notation we introduced in Section~\ref{sec:intro}, we define
\begin{equation}\label{potencial}
V_{\mu}(x)\!:=\int_1^{x+1}(u^p-u^q)du.
\end{equation}
\teoc{thm:familia} illustrates the application of the criticality results we have obtained in the previous section. This Theorem collects some of the cases that the results in~\cite{ManRojVil2015,ManRojVil2016-JDDE} did not cover. To prove the result we first need to show a technical lemma concerning the function $f(p)$ we have introduced in the introductory section. This Lemma ensures the uniqueness of the point $p_1$ in \teoc{thm:familia}. 

\begin{lema}\label{lema:tecnic}
The function $f(p)=(\tfrac{3+3p}{2})^{\frac{3+3p}{1+3p}}-2-3p$ has a unique zero on $(-\frac{1}{3},+\infty)$.
\end{lema}

\begin{prova}
We point out that the result is equivalent to show that $g(x)\!:=(\tfrac{x}{2})^{\frac{x}{x-2}}-x+1$ is monotonous decreasing on $(2,+\infty).$ Indeed, we have $g(x)=f(x/3-1)$, $\lim_{x\to 2}g(x)=e-1>0$ and $\lim_{x\to +\infty}g(x)=-\infty.$ The result follows then by continuity of the function $f$. Let us show the monotonicity of $g$. Elementary computations yield to
\[
g'(x)=-1+\left(\frac{x}{2}\right)^{\frac{x}{x-2}}\frac{x-2+\log(4)-2\log(x)}{(x-2)^2}
\]
and
\[
g''(x)=-\left(\frac{x}{2}\right)^{\frac{2}{x-2}}\frac{4+x^2-2x(2+\log(2)^2)+x\log(x)(\log(16)-2\log(x))}{(x-2)^4}.
\]
Since $\lim_{x\to 2^+}g'(x)<0$ it is enough to show that $g''<0$ in $(2,+\infty).$ On account of the expression of $g''$, this fact is equivalent to say that the function $\kappa (x)\!:=4+x^2-2x(2+\log(2)^2)+x\log(x)(\log(16)-2\log(x))$ is positive. This last statement is clear since one can easily verify that $\kappa(2)=\kappa'(2)=\kappa''(2)=0$ and $\kappa'''(x)=\frac{4\log(x/2)}{x^2}>0$ for all $x>2$. This shows the validity of the Lemma.
\end{prova}

\begin{prooftext}{Proof of \teoc{thm:familia}.}

Before start with the proof itself, notice that, since $p$ and $q$ are both different from $-1$, then the expression in~\eqref{potencial} writes
\begin{equation}\label{eq:Vx}
V_{\mu}(x)=\frac{(x+1)^{p+1}}{p+1}-\frac{(x+1)^{q+1}}{q+1}+h_0(\mu),
\end{equation}
where $h_0(\mu)\!:=\frac{p-q}{(p+1)(q+1)}$ corresponds to the energy level at the outer boundary of the period annulus of the family~\eqref{Xmu2} when $q>-1$. Observe that this is the situation in Theorem~\ref{thm:familia}. The projection of the period annulus on the $x$-axis is $\PI_{\mu}=(x_{\ell}(\mu),x_r(\mu))$, with
\[
x_{\ell}(\mu)=-1\text{ and }x_r(\mu)=\left(\frac{p+1}{q+1}\right)^{\frac{1}{p-q}}-1.
\]
Following the notation in \teoc{thm:criticality_finite}, from \refc{eq:Vx} and due to $p>q$, one can easily check that the family $\{h_0(\mu)-V_{\mu}(x)\}_{\mu\in\Lambda}$ is continuously quantifiable in any $\hat\mu=(\hat q,\hat p)\in\Lambda$ at $x=x_{\ell}$ by $\beta_{\ell}(\mu)$ with limit $b_{\ell}$ and at $x=x_r$ by $\beta_r(\mu)$ with limit $b_r$, where
\begin{equation}\label{eq:info_V}
\beta_{\ell}(\mu)=-(q+1),\ \beta_r(\mu)=-1,\ b_{\ell}=\frac{1}{\hat q+1}\ \text{and}\ b_r=V_{\hat\mu}'(x_r).
\end{equation}

Let us prove assertion in $(a)$ by applying \teoc{thm:criticality_finite} with $n=0$. To do so, let $\hat\mu=(\hat q,\hat p)$ with $\hat q=-\frac{1}{3}$ and $\hat p\in (-\frac{1}{3},+\infty)\setminus\{p_1\},$ where $p_1$ is the unique zero of $f(p)=(\tfrac{3+3p}{2})^{\frac{3+3p}{1+3p}}-2-3p$ on $(-\frac{1}{3},+\infty)$ (see Lemma~\ref{lema:tecnic}). In order to obtain the quantifier $\xi$ of $\bigl\{\PP[z\sqrt{h_0(\mu)}(g^{-1}_{\mu})''(z\sqrt{h_0(\mu)})]\bigr\}_{\mu\in\Lambda}$ we shall use the second part of~\cite[Theorem B]{ManRojVil2016-JDDE}. This result, for $n=0$, states that
$\xi(\mu)=-\min\{\bigl(\tfrac{\alpha_{\ell}}{\beta_{\ell}}\bigr)(\mu),\bigl(\tfrac{\alpha_{r}}{\beta_{r}}\bigr)(\mu)\}+1,$
where $\beta_{\ell}$ and $\beta_r$ are the functions in~\eqref{eq:info_V}, and $\alpha_{\ell}$ and $\alpha_r$ are the quantifiers of $\bigl\{(h_0-V_{\mu}){V_{\mu}}^{\frac{1}{2}}\mathscr R_{\mu}\bigr\}_{\mu\in\Lambda}$, with $\mathscr R_{\mu}\!:=\frac{(V_{\mu}')^2-2V_{\mu}V_{\mu}''}{(V_{\mu}')^3},$ at $x_{\ell}$ and $x_r$, respectively. Using the expression in~\eqref{eq:Vx}, it is a computation to show that the previous family is continuously quantifiable in $\hat\mu$ at $x=x_{\ell}$ by $\alpha_{\ell}(\mu)=q$ with limit $a_{\ell}=2\bigl(2-\frac{1}{\hat p+1}\bigr)^{\frac{3}{2}}$.
On the other hand, since $V_{\mu}$ is analytic at $x=x_r$ and $V_{\hat\mu}'(x_r)\neq 0,$ the family is continuously quantifiable in $\hat\mu$ at $x=x_r$ by $\alpha_r(\mu)=-1$ with limit 
\[
 a_r=\sqrt{h_0(\hat\mu)}\frac{V_{\hat\mu}'(x_r)^2-2h_0(\hat\mu)V_{\hat\mu}''(x_r)}{V_{\hat\mu}'(x_r)^{2}},
\]  
provided that $V_{\hat\mu}'(x_r)^2-2h_0(\hat\mu)V_{\hat\mu}''(x_r)\neq 0.$ This inequality is equivalent to require that $f(\hat p)\neq 0,$ which is satisfied since $\hat p\neq p_1$. Therefore, we have that $\bigl(\frac{\alpha_{\ell}}{\beta_{\ell}}\bigr)(\hat\mu)=\frac{1}{2}$ and $\bigl(\frac{\alpha_{r}}{\beta_{r}}\bigr)(\hat\mu)=1$ and so the quantifier of the family
$\bigl\{\PP[z\sqrt{h_0(\mu)}(g^{-1}_{\mu})''(z\sqrt{h_0(\mu)})]\bigr\}_{\mu\in\Lambda}$ at $z=1$ is $\xi(\mu)=\max\{\frac{q}{q+1},-1\}+1$.

We are in position to apply now \teoc{thm:criticality_finite}. To this end note that $\xi(\hat\mu)=\frac{1}{2}$. Then, the application of $(a)$ in \teoc{thm:criticality_finite} with $n=0$ shows that $\mathrm{Crit}\bigl((\out_{\hat\mu},X_{\hat\mu}),X_{\mu}\bigr)=0$. This proves the validity of $(a)$.

Let us now prove the assertion in $(b)$. In this case we consider $\hat\mu=(\hat q,\hat p)$ with $\hat q=0$ and $\hat p\in(0,+\infty)\setminus\{1\}$. By~\cite[Theorem E]{ManRojVil2015}, $\mathrm{Crit}\bigl((\Pi_{\hat\mu},X_{\hat\mu}),X_{\mu}\bigr)\geq 1$ and by~\cite[Theorem C]{ManRojVil2016-JDDE}, $\mathrm{Crit}\bigl((\Pi_{\hat\mu},X_{\hat\mu}),X_{\mu}\bigr)=1$ if $\hat p\in(0,+\infty)\setminus\{\tfrac{1}{2},1\}$. Then the result will follow by applying Theorem~\ref{thm:criticality_finite} with $n=1$ and $\hat\mu=(0,\tfrac{1}{2})$. From the computations in the proof of~\cite[Theorem C]{ManRojVil2016-JDDE} we have that the family 
\[
(\DD_{\boldsymbol\nu_1(\mu)}\circ\mathcal P)\bigl[z\sqrt{h_0(\mu)}(g_{\mu}^{-1})''(z\sqrt{h_0(\mu)})\bigr]
\]
with $\nu_1(\mu)=\frac{2q}{q+1}$ is continuously quantifiable in $\hat\mu=(0,\hat p)$ at $z=1$ by  $\xi(\hat\mu)=-\min\{\hat p-1,0\}=1-\hat p.$ In the case $\hat p=\tfrac{1}{2}$ we have $\xi(\hat\mu)=\frac{1}{2}$.
Then, applying \teoc{thm:criticality_finite}, we can assert that 
$\mathrm{Crit}((\Pi_{\hat\mu},X_{\hat\mu}),X_{\mu})\leq 1$ as desired.
\end{prooftext}

\subsection{Proof of Theorem~\ref{thm:loud}}
\begin{figure}
\centering
\includegraphics[scale=2]{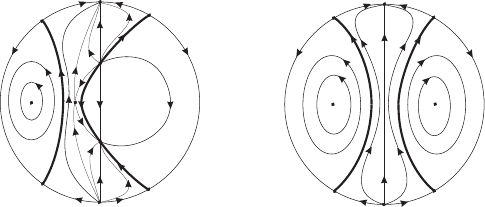}
\caption{\label{fig:retrat}Phase portrait of~\eqref{loud} in the Poincar\'e disc for $\mu=(D,F)\in\Lambda$ with $D<-1$ (left) and $D>-1$ (right), where for convenience we place the center at $(0,0)$ on the left of the centred invariant line $\{x=1\}$. The polycycle $\Pi_{\mu}$ at the outer boundary of the period annulus is the same in both cases: two hyperbolic saddles at infinity and the heteroclinic orbits between them. The invariant hyperbola $\mathcal C_{\mu}$ is in boldface type.}
\end{figure}

In this section we apply the techniques developed in Section~\ref{sec:criticality} to contribute in the study of the bifurcation of critical periodic orbits of the family of dehomogenized Loud's centers. We consider equation~\eqref{loud} with $\mu=(D,F)$ inside the set
\[
\Lambda\!:=\{(D,F)\in\R^2: F>1, D<0, D+F>0\}.
\]
It is known (see~\cite{MMV2} for instance) that if $F\notin\{0,1,\tfrac{1}{2}\}$ then the function
\[
H_{\mu}(x,y)=(1-x)^{-2F}\bigl(\tfrac{1}{2}y^2-q_{\mu}(x)\bigr)
\]
is a first integral of system~\eqref{loud}, where $q_{\mu}(x)=a(\mu)x^2+b(\mu)x+c(\mu)$ with
\[
a(\mu)=\frac{D}{2(1-F)},\ b(\mu)=\frac{D-F+1}{(1-F)(1-2F)} \text{ and } c(\mu)=\frac{F-D-1}{2F(1-F)(1-2F)}.
\]
The line at infinity, the conic $\mathcal C_{\mu}=\{\tfrac{1}{2}y^2-q_{\mu}(x)=0\}$ and the line $\{x=1\}$ are invariant curves of the system. If $\mu\in\Lambda$ then $\mathcal C_{\mu}$ is a hyperbola that intersects the $x$-axis at
\[
x=p_1(\mu)\!:=\frac{-b(\mu)-\sqrt{b(\mu)^2-4a(\mu)c(\mu)}}{2a(\mu)} \text{ and } x=p_2(\mu)\!:=\frac{-b(\mu)+\sqrt{b(\mu)^2-4a(\mu)c(\mu)}}{2a(\mu)},
\]
with $0<p_1(\mu)<p_2(\mu)$. For these parameters, the outer boundary of the period annulus of the center at the origin is formed by the union of the branch of the hyperbola $\mathcal C_{\mu}$ passing through the point $(p_1(\mu),0)$ and the line at infinity joining two hyperbolic saddles (see Figure~\ref{fig:retrat}). 

We are concerned with the criticality at the outer boundary of the period function for parameters $\hat\mu=(\hat D, 2)\in\Lambda$. As we advanced in the Introduction, by applying \cite[Lemma~14]{GGV} it follows that the change of coordinates
\[
(u,v)=(\phi(1-x),(1-x)^{-F}y), \text{ with }\phi(z)\!:=\frac{z^{-F}-1}{F},
\]
transforms the differential system~\eqref{loud} into the potential system
\begin{equation}\label{loud2}
\begin{cases}
\dot u = -v,\\
\dot v = (Fu+1)\bigl((Fu+1)^{-\frac{1}{F}}-1\bigr)\bigl(D(Fu+1)^{-\frac{1}{F}}-D-1\bigr).
\end{cases}
\end{equation}
This potential system has a non-degenerated center at the origin and the projection of its period annulus is the interval $\mathcal I_{\mu}\!:=(-\frac{1}{F},u_r(\mu))$, where
\[
u_r(\mu)\!:=\frac{(1-p_1(\mu))^{-F}-1}{F}.
\]
We point out that $\mathcal I_{\mu}$ is the image by $x\mapsto \phi(1-x)$ of $(-\infty,p_1(\mu))$. Setting $z=\phi^{-1}(u)=(Fu+1)^{-1/F}$ one can check that
\begin{equation}\label{formulas}
 \begin{array}{ll}
  V_{\mu}(u)=h_0(\mu)-z^{-2F}V_0(z,\mu) & \text{ with }V_0(z,\mu)=\frac{D}{2-2F}z^2+\frac{1+2D}{2F-1}z-\frac{D+1}{2F}, \\[3pt]
  V_{\mu}'(u)=z^{-F}V_1(z,\mu) & \text{ with }V_1(z,\mu)=(z-1)(D(z-1)-1),  \\[3pt]
  V_{\mu}''(u)=V_2(z,\mu) & \text{ with }V_2(z,\mu)=D(F-2)z^2-(2D+1)(F-1)z+F(D+1),  \\[3pt]
  V_{\mu}^{(3)}(u)=z^FV_3(z,\mu)& \text{ with }V_3(z,\mu)=-2D(F-2)z^2+(2D+1)(F-1)z,  \\[3pt]
   V_{\mu}^{(4)}(u)=z^{2F}V_4(z,\mu)& \text{ with }V_4(z,\mu)=2D(F^2-4)z^2-(2D+1)(F^2-1)z,
 \end{array}
\end{equation}
where $h_0(\mu)\!:=\frac{F-D-1}{2F(F-1)(2F-1)}$ is the energy level at the outer boundary of the period annulus for all $\mu\in\Lambda$.

Before the proof of Theorem~\ref{thm:loud} a required technical result is shown. Its proof follows a similar argument than~\cite[Lemma~3.4]{RojVil2018}.

\begin{lema}\label{tecnic2}
Let $\mu=(D,2)\in\Lambda$. If $D\neq-\frac{1}{2}$ then $V_{\mu}'(u_r(\mu))^2+\frac{D-1}{6}V_{\mu}''(u_r(\mu))\neq 0$. 
\end{lema}

\begin{prova}
Let us denote $L(z,\mu)\!:=z^{-2F}V_1(z,\mu)^2+\frac{D-1}{6}V_2(z,\mu)$. From the expressions in~\eqref{formulas} we have
\[
V_{\mu}'(u)^2+\frac{D-1}{6}V_{\mu}''(u)=L(\phi^{-1}(u)).
\]
Then the result will follows once we show that $L(z,\mu)$ at $z=\phi^{-1}(u_r(\mu))=1-p_1(\mu)=1+\frac{b+\sqrt{b^2-4ac}}{2a}$ does not vanish for $D\neq -\frac{1}{2}$.

We claim that if $\hat\mu=(\hat D,2)\in\Lambda$ is a zero of $L(1-p_1(\mu),\mu)$ then the derivative of the map $D\mapsto L(1-p_1(\mu),\mu)$ is strictly negative at $\mu=\hat\mu$. Indeed, first we notice that if $\hat\mu$ is a zero of $L(1-p_1(\mu),\mu)$ then
\begin{equation}\label{zeroL}
(1-p_1(\hat\mu))^{-2F}=\frac{(1-D)V_2(1-p_1(\hat\mu),\hat\mu)}{6V_1(1-p_1(\hat\mu),\hat\mu)}.
\end{equation}
In addition,
\begin{align*}
\frac{d}{dD}L(1-p_1(\mu),\mu)&=z^{-2F}\left.\left(2V_1(z,\mu)\partial_z V_1(z,\mu)-\frac{2FV_1(z,\mu)^2}{z}\right)\right|_{z=1-p_1(\mu)}\partial_D(1-p_1(\mu))\\
&\phantom{=}+z^{-2F}\partial_D(V_1(z,\mu)^2)+\frac{1}{6}\left.\partial_D((D-1)V_2(z,\mu))\right|_{z=1-p_1(\mu)}.
\end{align*}
Substituting the equality~\eqref{zeroL} in the previous identity and evaluating at $\hat\mu=(\hat D,2)$, we obtain that $\frac{d}{dD}L(1-p_1(\hat\mu),\hat\mu)$ can be written as an algebraic expression in $\hat D$,
\[
\frac{d}{dD}L(1-p_1(\hat\mu),\hat\mu)=\frac{r_1(\hat\mu)\sqrt{\Delta(\hat\mu)}+r_2(\hat\mu)}{r_3(\hat\mu)\sqrt{\Delta(\hat\mu)}+r_4(\hat\mu)} \text{ with } \Delta(\mu)\!:=(b^2-4ac)(\mu)
\]
and $r_i\in\R[\mu]$ satisfying
\[
r_1(\hat\mu)^2\Delta(\hat\mu)-r_2(\hat\mu)^2=-\frac{2187}{2}\hat D^4(\hat D+1)^3(\hat D+2)(108-486\hat D+810\hat D^2-594\hat D^3+162\hat D^4)
\]
and
\[
r_3(\hat\mu)^2\Delta(\hat\mu)-r_4(\hat\mu)^2=-5668704\hat D^6(\hat D-1)^2(\hat D+1)^3(\hat D+2)^3.
\]
We point out that if $\hat D\in(-2,0)\setminus\{-1\}$ the previous two expressions does not vanish. Indeed, the first one follows by Descartes' rule of signs while the second is straightforward. This proves that $\frac{d}{dD}L(1-p_1(\hat\mu),\hat\mu)$ is well defined and non-vanishing for all $\hat\mu=(\hat D,2)\in\Lambda$ with $\hat D\neq -1$. In addition, one may explicitly check that
\[
\frac{d}{dD}L(1-p_1(\hat\mu),\hat\mu)=-\frac{1}{12} \text{ if }\hat\mu=(-1,2).
\]
This shows that $\frac{d}{dD}L(1-p_1(\hat\mu),\hat\mu)$ is negative for all $\hat D\in(-2,0)$, proving the validity of the claim.

The claim implies that the map $D\mapsto L(1-p_1(\mu),\mu)$ as at most one zero for $D\in(-2,0)$ and $F=2$. It is a simple computation to show that $D=-\frac{1}{2}$ is a zero and this proves the result.
\end{prova}

\begin{prooftext}{Proof of Theorem~\ref{thm:loud}.}
The strategy of the proof will be to use Theorem~\ref{thm:criticality_finite} with $n=1$ and $\hat\mu=(\hat D,2)$ with $\hat D\in(-2,0)\setminus\{-1/2\}$. As in the previous section, in order to obtain the quantifier $\xi$ of $\{(\DD_{\boldsymbol\nu_1(\mu)}\circ\mathcal P)\bigl[z\sqrt{h_0(\mu)}(g_{\mu}^{-1})''(z\sqrt{h_0(\mu)})\bigr]\}_{\mu\in\Lambda}$ at $z=1$ in $\hat\mu$ we shall use the second part of~\cite[Theorem~B]{ManRojVil2016-JDDE}. In this case, for $n=1$, it states that 
\begin{equation}\label{xi-loud}
\xi(\mu)=-\min\left\{\left(\frac{\alpha_{\ell}}{\beta_{\ell}}\right)(\mu),\left(\frac{\alpha_{r}}{\beta_{r}}\right)(\mu)\right\}-\frac{1}{2}\nu_1(\mu),
\end{equation}
where $\alpha_{\ell}$ and $\alpha_r$ are the quantifiers of
\[
\Psi_{\mu}(u)\!:=\frac{1}{V_{\mu}'(u)} W\left[\left(\frac{V_{\mu}}{h_0(\mu)-V_{\mu}}\right)^{\frac{1}{2}\nu(\mu)},(h_0(\mu)-V_{\mu})V_{\mu}^{\frac{1}{2}}\mathcal R_{\mu} \right]\!(u)
\]
with $\mathcal R_{\mu}\!:= \frac{(V_{\mu}')^2-2V_{\mu}V_{\mu}''}{(V_{\mu}')^3}$, and $\beta_{\ell}$ and $\beta_r$ are the quantifiers of $(h_0-V_{\mu})$, at $u=-1/F$ and $u=u_{r}(\mu)$, respectively. The map $\nu:\Lambda\rightarrow\R$ is a continuous function to be determined in such a way that $\Psi_{\mu}$ is continuously quantifiable at $u=-1/F$ for $F\approx 2$.

Firstly, in \cite[Lemma~3.3]{RojVil2018} it was shown that 
\begin{equation}\label{quant1}
h_0(\mu)-V_{\mu}(u)\sim_{-\frac{1}{F}} \frac{D}{2-2F}(Fu+1)^{\frac{2F-2}{F}} \text{ at }\hat\mu
\end{equation}
and
\[
h_0(\mu)-V_{\mu}(u)\sim_{u_r(\mu)} V_{\mu}'(u_r(\mu))(u_r(\mu)-u)\text{ at }\hat\mu
\]
with $V_{\mu}'(u_r(\mu))\neq 0$. With the notation introduced before, we have then 
\[
\beta_{\ell}=\frac{2-2F}{F} \text{ and }\beta_r=-1.
\]
Secondly, a computation shows that
\begin{equation}\label{Psi}
\Psi_{\mu}(u)=\frac{1}{2\sqrt{V_{\mu}(u)}V_{\mu}'(u)^5}\left( \frac{V_{\mu}(u)}{h_0(\mu)-V_{\mu}(u)}\right)^{\frac{\nu(\mu)}{2}}f_{\mu}(u),
\end{equation}
where $f_{\mu}\!:=4V_{\mu}^2(V_{\mu}-h_0(\mu))V_{\mu}'V_{\mu}'''-((V_{\mu}')^2-2V_{\mu}V_{\mu}'')((V_{\mu}')^2(3V_{\mu}+h_0(\mu)(\nu(\mu)-1))-6(V_{\mu}-h_0(\mu))V_{\mu}V_{\mu}'')$.

To compute $\alpha_{\ell}$, taking~\eqref{formulas} into account, one can verify with the help of an algebraic manipulator that
\[
f_{\mu}(u)=\left.\frac{z^{-6F}}{2(F-1)^3F^3(2F-1)^3}\psi_{\mu}(z)\right|_{z=\phi^{-1}(u)},
\]
where $\psi_{\mu}$ is the sum of $25$ monomials of the form $\kappa(\mu)z^{n_1+n_2F}$ with $n_i\in\Z$ for $i=1,2$, and $\kappa:\Lambda\rightarrow\R$ a well-defined rational function at $\hat\mu$. Moreover, the monomial with highest exponent for $\mu\approx\hat\mu$ is $D^3F(1+D-F)^2(F-2)(2F-1)(2+F(\nu(\mu)-1)-\nu(\mu))z^{6+4F}$. We point out that the coefficient of the previous monomial vanishes at $F=2$ so $\psi_{\mu}$ is not continuously quantifiable at $z=\infty$ in $\hat\mu$ for arbitrary function $\nu(\mu)$. In order to succeed with the continuous quantification, we must take
\[
\nu(\mu)\!:=\frac{F-2}{F-1}.
\]
In this case, the previous coefficient vanishes for all $\mu\approx\hat\mu$ and the monomial with highest exponent is $D^2F(1+2D)(1+D-F)^2(F^2-2F-3)z^{5+4F}$. Consequently,
\[
f_{\mu}(u)\sim_{-\frac{1}{F}} \frac{D^2(1+2D)(1+D-F)^2(F^2-2F-3)}{2F^2(F-1)^3(2F-1)}(Fu+1)^{\frac{2F-5}{F}} \text{ at }\hat\mu.
\]
In addition, using the expressions in~\eqref{formulas}, we have that $V_{\mu}'(u)\sim_{-\frac{1}{F}} D(Fu+1)^{\frac{F-2}{F}}$. Using the previous quantifications together with~\eqref{quant1} in the expression~\eqref{Psi} we have
\[
\Psi_{\mu}(u)\sim_{-\frac{1}{F}} C_1(\mu)(Fu+1)^{\frac{7-4F}{F}}\text{ at }\hat\mu,
\]
where
\[
C_1(\mu)\!:=\frac{(1+2D)(1+D-F)^2(F^2-2F-3)}{4D^3F^2(F-1)^3(2F-1)\sqrt{h_0(\mu)}}\left(\frac{(2-2F)h_0(\mu)}{D}\right)^{\frac{F-2}{2(F-1)}}.
\]
Finally, by simple computations from~\eqref{Psi} using that $u=u_r(\mu)$ is a regular value of $V_{\mu}$ and a simple zero of $h_0(\mu)-V_{\mu}(u)$ with $V_{\mu}'(u_r)\neq 0$, we have that
\[
\Psi_{\mu}(u)\sim_{u_r(\mu)} C_2(\mu) (u_r(\mu)-u)^{\frac{F-2}{2(F-1)}}\text{ at }\hat\mu,
\] 
where
\[
C_2(\mu)\!:=\frac{(D-1)(6V_{\mu}'(u_r(\mu))^2+(D-1)V_{\mu}''(u_r(\mu)))}{12\sqrt{3(1-D)}V_{\mu}'(u_r(\mu))}.
\]
We point out that $\nu(\mu)=\frac{F-2}{F-1}$ and that $6V_{\hat\mu}'(u_r(\hat\mu))^2+(D-1)V_{\hat\mu}''(u_r(\hat\mu))$ does not vanish if $\hat D\neq -1/2$ by Lemma~\ref{tecnic2}. In particular,
\[
\alpha_{\ell}=\frac{4F-7}{F} \text{ and }\alpha_r=\frac{2-F}{2(F-1)}.
\]

Consequently, using the equality in~\eqref{xi-loud}, we have that $\{(\DD_{\boldsymbol\nu(\mu)}\circ\mathcal P)\bigl[z\sqrt{h_0(\mu)}(g_{\mu}^{-1})''(z\sqrt{h_0(\mu)})\bigr]\}_{\mu\in\Lambda}$ with $\nu(\mu)=\frac{F-2}{F-1}$ is continuously quantifiable at $z=1$ in $\hat\mu$ by 
\[
\xi(\mu)=\left.-\min\left\{\frac{4F-7}{2-2F},\frac{F-2}{2(F-1)}\right\}-\frac{F-2}{2(F-1)}\right|_{F=2}=\frac{1}{2}.
\]
The result follows then applying Theorem~\ref{thm:criticality_finite} with $n=1$, proving that $\mathrm{Crit}((\Pi_{\hat\mu},X_{\hat\mu}),X_{\mu})\leq 1$ as desired.
\end{prooftext}

%\section{Offtopic}
%[ESTE LEMA NO SE USA, PERO ES INTERESANTE Y QUIZAS UTIL DESPUES]
%
%\begin{lema}\label{lemma:M1}
%Let $f_{\alpha}(x)\sim_{\infty} ax^{\alpha}$ in a neighborhood of $\alpha=-1$ and assume that $M_1[f_{\alpha}]\neq 0$ if $\alpha<-1$. Then 
%\[
%\lim_{\alpha\rightarrow -1^+}(\alpha+1) M_1[f_{\alpha}]=-a.
%\]
%\end{lema}
%\begin{prova}
%Let us fix $\epsilon>0$. Since $f_{\alpha}(x)\sim_{\infty} ax^{\alpha}$, there exist $M>0$ and $\alpha_1<-1$ such that $\abs{ x^{-\alpha}f_{\alpha}(x)-a}< \frac{\epsilon}{8}$ for all $x>M$ and $\alpha\in (\alpha_1,-1)$. Let us take $N\!:=\max\{ \abs{f_{\alpha}(x)} : x\in[0,M],\alpha\in [\alpha_1,-1]\}$. 
%On the one hand, there exists $\alpha_2\in(\alpha_1,-1)$ such that 
%\[
%\abs{(\alpha+1)\int_0^{M}f_{\alpha}(x)dx} \leq \abs{\alpha+1}MN\leq \frac{\epsilon}{2} \text{ for all }\alpha\in(\alpha_2,-1).
%\]
%On the other hand, there exists $\alpha_0\in(\alpha_2,-1)$ such that $M^{1+\alpha}<2$ and $\abs{1-M^{1+\alpha}}<\frac{\epsilon}{4\abs{a}}$ for all $\alpha\in(\alpha_0,-1)$. Therefore,
%\begin{align*}
%\abs{(\alpha+1)\int_M^{\infty}f_{\alpha}(x)dx+a}
%&\leq 
%\abs{\alpha+1}\int_M^{\infty}\abs{x^{-\alpha}f_{\alpha}(x)-a}x^{\alpha}dx + \abs{(\alpha+1)\int_M^{\infty} ax^{\alpha}dx +a}\\
%&\leq 
%\frac{\epsilon}{8} M^{\alpha+1}+\abs{a}\abs{1- M^{\alpha+1}}\leq \frac{\epsilon}{2}
%\end{align*}
%for all $\alpha\in(\alpha_0,-1)$. Consequently,
%\[
%\abs{(\alpha+1)M_1[f_{\alpha}]+a}\leq \abs{(\alpha+1)\int_0^{M}f_{\alpha}(x)dx} + \abs{(\alpha+1)\int_M^{\infty}f_{\alpha}(x)dx+a}\leq \epsilon
%\]
%for all $\alpha\in(\alpha_0,-1)$ as we desired.
%\end{prova}
%
\section*{Acknowledgements}

The author want to thank Prof. Pavao Marde{\v{s}}i{\'c} for the short stay of four months at the Institut de Math\'ematiques de Bourgogne, Dijon, that gave rise to the start of this work; and to Prof. Francesc Ma\~nosas and Dr. Jordi Villadelprat for the fruitful discussions that led to some corrections of this paper. The author is partially supported by the
MINECO/FEDER grants MTM2017-82348-C2-1-P and MTM2017-86795-C3-1-P.

%\bibliography{mybibfile}

\end{document}